% Version 3/8/00

%\magnification \magstep1
\input amstex
\documentstyle{amsppt}
\topmatter
\rightheadtext{Finite jet determination}
\leftheadtext{Peter Ebenfelt}
\date{\number\month-\number\day-\number\year}\enddate

%\loadeufm
%\define \p{\eufm p}
%\define \q{\eufm q}
%\define \m{\eufm m}
%\define \ass{\text{\rm Ass }}
%\define \sol{\text{\rm Sol}}

\define \bR{\Bbb R}
%\define \bR{{\bold R}}

%\define \bZ{{\bold Z}}
\define \bC{\Bbb C}
%\define \bC{{\bold C}}

%\define \ssneq{\subsetneqq}

\define \hol{\text{\rm hol}}

\def\dim {\text {\rm dim}}

\def\aut{\text{\rm aut}}
\def\Aut{\text{\rm Aut}}
\def\pb11{\overline{p_{11}}}
\def\pb02{\overline{p_{02}}}

\title Finite jet determination of holomorphic mappings at the boundary\endtitle
\author Peter Ebenfelt\footnote{Supported in part by a grant from the
Swedish Natural Science Research
Council.\newline}
\endauthor
%\affil Department of Mathematics\\Royal Institute of Technology\\100
%44 Stockholm, Sweden\endaffil
\address Department of Mathematics, Royal Institute of Technology, 100
44 Stockholm, Sweden\endaddress
\email ebenfelt\@math.kth.se\endemail

\abstract Let $M,M'\subset\bC^N$ be smooth real hypersurfaces and
assume that $M$ is $k_0$-nondegenerate at $p_0\in M$. We prove that
holomorphic mappings that extend smoothly to $M$, sending a
neighborhood of $p_0$ in $M$ diffeomorphically into $M'$ are
completely determined by (and depend smoothly on) their $2k_0$-jets at $p_0$. As an application
of this result, we give sufficient conditions on a smooth real
hypersurface which guarantee that the space of
infinitesimal CR automorphisms is finite dimensional.
\endabstract


\endtopmatter
\document

\heading 0. Introduction\endheading

A classical theorem of H. Cartan ([HCa]) states that an automorphism
$f$ of a bounded domain $D\subset\bC^N$ is completely determined
by its $1$-jet, i.e. its value and derivatives of order one, at
any point $Z_0\in D$. If $D$, in addition, is assumed to be
smoothly ($C^\infty$) bounded and strictly pseudoconvex, then by
Fefferman's theorem [Fe] any such automorphism extends smoothly
to the boundary $\partial D$ as an automorphism $\partial D\to
\partial D$. It is then natural to ask: is $f$ completely
determined by a finite jet at a boundary point $p\in
\partial D$? An affirmative answer to this question, when $D$ is
strictly pseudoconvex, follows from the
work of Chern and Moser [CM] (see also E. Cartan [ECa1--2] for the
case $N=2$, and Tanaka [T1--2]). Indeed, the following local version
of Cartan's theorem is a consequence of their work. {\it Any
holomorphic mapping which is defined locally on one side of a
smooth, Levi nondegenerate real hypersurface $M\subset \bC^N$ and
extends smoothly to $M$, sending $M$ diffeomorphically into
another smooth real hypersurface $M'\subset \bC^N $, is completely
determined by its $2$-jet at a point $p\in M$.} Observe that the
conclusion is nontrivial even in the strictly pseudoconcave case
when the mapping extends holomorphically to a full neighborhood of
$p$.

The main objective of the present paper is to extend the above
mentioned local result to a more general class of real
hypersurfaces (Theorem 1 below). We should point out that the
result for Levi nondegenerate hypersurfaces follows from the
construction of a unique Cartan connection on a certain principal
$G$-bundle over such a hypersurface. There is no analogue of this
construction in the more general situation considered in this paper.

Let $M\subset\bC^N$ be a smooth real hypersurface
and assume that $M$ is defined locally near a point $p_0\in M$ by
the equation $\rho(z,\bar z)=0$, where $\rho$ is a smooth function
with $\rho(p_0,\bar p_0)=0$ and $d\rho(p_0,\bar p_0)\neq 0$. Let
$L_{\bar 1},\ldots, L_{\bar n}$, with $n=N-1$, be a basis for the CR
vector fields on $M$.  We shall say that $M$ is {\it
$k_0$-nondegenerate} at $p_0$ if $$\text{\rm Span} \left\{
L^{\bar \alpha}\rho_z(p_0,\bar
p_0)\:|\alpha|\leq k_0\right \}=\bC^N,\tag 0.1$$ where
$\rho_z:=(\partial_{z_j} \rho)_{1\leq j\leq N}$, $\partial_
{z_j}:=\partial/\partial {z_j}$, and
standard multi-index notation for differential operators is used
i.e.\ $L^{\alpha}:=L_{ \bar 1}^{\bar \alpha_1}\ldots L_{\bar
n}^{\bar \alpha_n}$. This nondegeneracy condition will be given in
a different, but equivalent, form in terms of the intrinsic
geometry of $M$ in the next section. The reader is referred to the
book [BER3] for basic material on real submanifolds in complex
space and CR structures, and further discussion of various
nondegeneracy conditions (see also \S1 of the present paper). We
mention here only that Levi
nondegeneracy at a point $p_0\in M$ is equivalent to
$1$-nondegeneracy. Our main result is the following.

\proclaim{Theorem 1} Let $M, M'\subset\bC^{N}$ be smooth
($C^\infty$) real hypersurfaces. Let $f,g\:U\to \bC^N$, where
$U\subset \bC^N$ is an open connected subset with $M$ in its
boundary, be holomorphic mappings which extend smoothly to $M$ and
send $M$ diffeomorphically into $M'$. If $M$ is
$k_0$-nondegenerate at a point $p_0\in M$ and $$
(\partial_z^\alpha f)(p_0)=(\partial_z^\alpha g)(p_0),\quad
\forall\alpha\in\Bbb Z_+^N\: |\alpha|\leq 2k_0,\tag 0.2 $$ then $f\equiv g$ in $U$.
\endproclaim

Finite jet
determination of holomorphic mappings sending one real submanifold
into another has attracted much attention in recent years. We mention here
the papers [BER1--2, 4--5],
[L], [Han1--2], [Hay], [Z]. The reader is also referred to
the survey article [BER6] for a more detailed history. However, in all
the above mentioned papers, it is either assumed that $M$ and $M'$ are
real-analytic (which will imply that all mappings $f$ extend
holomorphically to some neighborhood of $M$), or the conclusion is
that the formal power series of the mapping $f$ is determined
by a finite jet (see [BER4], [L]). Theorem 1 appears to be, to the best of the
author's knowledge, the first finite determination result, since the work of
Chern and Moser mentioned above, which applies to merely smooth
hypersurfaces and smooth mappings. We should mention
that if $M$ and $M'$ are
real-analytic, then the conclusion of Theorem 1 was proved in
[BER2] (cf.\ also [Han1] and [Z]). A related notion is that of unique
continuation at the boundary for holomorphic mappings. A unique
continuation principle is said to hold for a class of mappings at a
point $p$ if any mapping from this class which agrees with the
constant mapping to infinite
order at $p$ is necessarily constant. (Observe that, due to the
nonlinear nature of mapping problems, a unique continuation principle
for a class of mappings into a manifold does not imply that two
mappings, in this class, which agree to infinite order are necessarily
the same.) We shall not address this problem further here. We
mention the papers [ABR], [BR], [BL], [CR], [E1], [HK], and refer
the interested reader to these papers for further information.

The proof of Theorem 1 is
based on
Theorem 2 below, and
a result from [BER4], alluded to above, which asserts that,
under the assumptions of Theorem 1, the jet of $f$ at $p_0$ of any order
is completely determined by its $2k_0$-jet. The proof of Theorem 1 is
given at the end of \S 3.

Our second result, which is the basis for Theorem 1 above, states,
loosely speaking, that given two suitably nondegenerate real
hypersurfaces, there is a system of differential equations, which
is complete in a certain sense, such that any CR diffeomorphism
$f\:M\to M'$ must satisfy this system. The idea to look for such a
differential system goes back to the work of E. Cartan and
Chern--Moser mentioned above. The approach was further developed
in the work of Han. To formulate the result more precisely, we
need to fix some notation. Let us denote by $J^k(M,M')_{(p,p')}$
the space of $k$-jets at $p\in M$ of smooth mappings $f\:M\to M'$
with $f(p)=p'\in M'$. Given coordinate systems $x=(x_1,\ldots,
x_{2N-1})$ and $x'=(x'_1,\ldots, x'_{2N-1})$ on $M$ and $M'$ near
$p$ and $p'$, respectively, there are natural coordinates
$\lambda^k:=(\lambda^\beta_i)$, where $1\leq i\leq 2N-1$ and
$\beta\in\Bbb Z_+^{2N-1}$ with $1\leq |\beta|\leq k$, on
$J^k(M,M')_{(p,p')}$ in which the $k$-jet at $p$ of a smooth
mapping $f\:M\to M'$ is given by
$\lambda_i^\beta=(\partial_x^\beta f_i)(p)$, $1\leq |\beta|\leq k$
and $1\leq i\leq 2N-1$.

\proclaim{Theorem 2} Let $M, M'\subset\bC^{N}$ be smooth
($C^\infty$) real hypersurfaces. Assume that $M$ is
$k_0$-nondegenerate at a point $p_0$. Let $f^0\:M\to M'$ be a
smooth CR diffeomorphism. Then, for any multi-index $\alpha\in\Bbb
Z_+^{2N-1}$ with $|\alpha|=k_0^3+k_0^2+k_0+2$ and any
$j=1,\ldots,2N-1$, there are smooth functions
$r^\alpha_j(\lambda^{k};x')(x)$ on $U$, where $k:=
k_0^3+k_0^2+k_0+1$ and $U\subset
J^{k}(M,M')_{(p_0,p_0')}\times M\times M'$ is an open
neighborhood of $((\partial_x^\beta f^0)(0),f^0(p_0),p_0)$, such
that $$
\partial_x^\alpha f_j=r^\alpha_j(\partial_x^\beta f;f),\quad \forall
|\alpha|= k_0^3+k_0^2+k_0+2,\ j=1\ldots, 2N-1,\tag 0.3 $$ where
$1\leq |\beta|\leq k$, for every smooth CR diffeomorphism
$f\:V\to M'$, where $V\subset M$ is some open neighborhood of $p_0$,
with $((\partial_x^\beta
f)(0),f(p_0),p_0)\in U$. Moreover, the functions $r^\alpha_j$ are
rational in $\lambda^{k}\in J^{k}(M,M')_{(p_0,p_0')}$; here,
$x=(x_1,\ldots,x_{2N-1})$ and $x'=(x'_1,\ldots, x'_{2N-1})$ are
any local coordinate systems on $M$ and $M'$ near $p_0$ and
$f^0(p_0)$, respectively, and $f_i:=f\circ x'_i$.
\endproclaim

Similar results for real-analytic hypersurfaces can also be
found in [Han1--2] and [Hay]. The idea behind the proof of Theorem 2
is to consider the tangent mapping $df\:\bC TM\to \bC TM'$ and
derive differential equations for $df$ using properties of a
sequence of invariant tensors (generalized Levi forms) which were
developed in the author's paper [E3]. The proof of Theorem 2 is given in \S 3.

We conclude this introduction by giving two applications of
Theorems 1 and 2. For this, we need some more notation. A smooth
real vector field $X$ on $M$ is called an {\it infinitesimal CR
automorphism} if the local 1-parameter group of diffeomorphisms,
$\exp tX$, generated by $X$ is a local group of CR diffeomorphisms
(see e.g.\ [BER2] or [S1--2]). The set of infinitesimal CR
automorphisms, defined near $p\in M$, forms a vector space over
$\bR$ denoted by $\aut(M,p)$. We shall give a sufficient condition
on $M$ at a point $p$ for $\dim_\bR\aut(M,p)<\infty$. A smooth
real hypersurface $M\subset \bC^N$ is called (formally) {\it
holomorphically degenerate} at $p\in M$, if there exists a formal
holomorphic vector field $$ Y=\sum_{j=1}^N
a_j(z)\partial_{z_j},\tag 0.4 $$ where the $a_j(z)$ are formal
power series in $z-p$, which is tangent to $M$, i.e.\ such that
the Taylor series at $p$ of a defining function $\rho(z,\bar z)$
for $M$ divides $(Y\rho)(z,\bar z)$ in the ring of formal power
series in $(z-p,\bar z-\bar p)$. Being holomorphically
nondegenerate (i.e.\ the opposite of being degenerate) at a point
is a strictly weaker condition than that of being
$k$-nondegenerate for some integer $k$. (See [BER3, Chapter XI]
for a more detailed description of the relationship between the
two notions). Also, recall that $M$ is said to be minimal at $p\in
M$ (in the sense of Tumanov and Trepreau) if $M$ does not contain
a complex hypersurface through $p$.

\proclaim{Theorem 3} Let $M\subset\bC^N$ be a smooth ($C^\infty$) real
hypersurface which is holomorphically nondegenerate and
minimal at $p_0$. Then,
$$
\dim_\bR\aut(M,p_0)\leq (2N-1)\binom{4N-3}{2N-2}\tag 0.5
$$
\endproclaim

A real-analytic hypersurface $M$ is said to be holomorphically
degenerate at $p\in M$ if there exists a holomorphic vector field,
i.e. a vector field of the form \thetag{0.4} with the $a_j(z)$
holomorphic, tangent to $M$ near $p$. This definition turns out to
be equivalent to the one given in the smooth category above (i.e.
using formal vector fields) for a real-analytic hypersurface (see
[BER3, Proposition 11.7.4]). Stanton [S2] proved that
$\dim_\bR\hol(M,p)<\infty$ for a real-analytic hypersurface $M$,
where $\hol(M,p)$ denotes the subspace of $\aut(M,p)$ consisting
of those infinitesimal CR automorphisms which are real-analytic,
if and only if $M$ is holomorphically nondegenerate at $p$. The
corresponding statement (as well as results for higher
codimensional real-analytic submanifolds) for $\aut(M,p)$, with
$M$ real-analytic, was proved in [BER2]. In contrast to the
real-analytic case, the condition of (formal) holomorphic
nondegeneracy is not necessary in Theorem 3. A real smooth
hypersurface $M$ in $\bC^2$ which is holomorphically degenerate
and minimal at $0$, but  everywhere Levi nondegenerate outside $0$
is given in [BER3, Example 11.7.29]. The fact that $M$ is Levi
nondegenerate outside $0$ can be seen to imply (see the concluding
remarks in \S4.2) that $\dim_\bR\aut(M,0)$ satisfies the bound in
\thetag{0.5}. However, if there exists a vector field $$
Y=\sum_{j=1}^N a_j(z,\bar z)\partial_{z_j},\tag 0.6 $$ where the
$a_j(z,\bar z)$ are smooth functions whose restrictions to $M$ are
CR, tangent to $M$ near $p$, then the arguments in [S2] easily
show that $\dim_\bR\aut(M,0)=\infty$. This discrepancy is
addressed further in \S4.2. The proof of Theorem 3 is given in
\S3.

For our final result, we shall denote by $\Aut(M,p)$ the stability
group of $M$ at $p\in M$, i.e. the group of germs at $p$ of local
CR diffeomorphisms $f\: V\to M$, where $V\subset M$ is some open
neighborhood of $p$, with $f(p)=p$. If $M$ is $k_0$-nondegenerate at
$p_0$, then, by Theorem 1, the jet mapping $j^{2k_0}_p$ sends
$\Aut(M,p_0)$ injectively into the jet group
$G^{2k_0}(\bC^N)_{p_0}\subset J^{2k_0}(\bC^N,\bC^N)_{(p_0,p_0)}$,
which 
consists of those jets that are invertible at $p_0$. 
We shall show that
the elements of $\Aut(M,p_0)$
depend smoothly on their $2k_0$-jets at $p_0$. More precisely, we have the
following result.

\proclaim{Theorem 4} Let $M\subset \bC^N$ be a smooth ($C^\infty$)
real hypersurface which is $k_0$-nondegenerate at $p_0\in M$.
Then, the jet mapping $$j^{2k_0}\:\Aut(M,p_0)\to
G^{2k_0}(\bC^N)_{p_0}$$ is injective and, for every
$f^0\in\Aut(M,p_0)$, there exist an open neighborhood $U_0$ of
$j^{2k_0}_{p_0}(f^0)$ in $G^{2k_0}(\bC^N)_{p_0}$, an open
neighborhood $V_0$ of $p_0$ in $M$, and a smooth ($C^\infty$)
mapping $F\:U_0\times V_0\to M$ such that
$$F(j^{2k_0}_{p_0}(f),\cdot)=f,\tag 0.7 $$ for every $f\in
\Aut(M,p_0)$ with $j^{2k_0}_{p_0}(f)\in U_0$.
\endproclaim

For real-analytic hypersurfaces, the result in Theorem 4 (with
real-analytic dependence) was proved in [BER1]. (See [BER4] for
the higher codimensional case; cf. also [Z].). 

\remark{Acknowledgement} The author would like to thank B. Lamel and D. Zaitsev for
many helpful comments and discussions on a preliminary version of this paper.
\endremark

\heading 1. Preliminaries \endheading

A real hypersurface $M\subset \bC^N$ inherits a CR structure $\Cal V:=
T^{0,1}\bC^N\cap \bC T M$ from
the ambient complex space $\bC^N$. (Here, $T^{0,1}\bC^N$ denotes the
usual bundle of
$(0,1)$ vectors in $\bC^N$.) In this section, we shall
consider abstract, not necessarily embedded (or integrable), CR
structures. At the end of this section, we shall again specialize
to embedded hypersurfaces, which substantially simplifies some of
the computations in subsequent sections. The reader is referred to the
concluding remarks in \S4 for a brief discussion of the abstract
case.

Let $M$ be a smooth ($C^\infty$) manifold with a CR structure $\Cal
V\subset\bC TM$. Recall that this means that $\Cal V$ is a formally
integrable subbundle (the commutator of two sections of $\Cal V$ is
again a section of $\Cal V$) such that $\Cal V_p\cap\bar \Cal V_p=\{0\}$ for
every $p\in M$. Sections of
the CR bundle are called {\it CR vector fields}. We shall denote by
$n\geq 1$ the CR dimension of the
CR manifold $M$, which by definition is the complex fiber dimension of
$\Cal V$, and we shall assume that the CR structure is of
hypersurface type, i.e. that $\dim_\bR M=2n+1$. The reader is referred
to [BER3] for an introduction to CR structures.

We define
two subbundles
$T^0M\subset T'M\subset \bC T^*M$ as
follows
$$
T^0M:=(\Cal V\oplus\bar\Cal V)^\perp,\quad T'M=\Cal V^\perp,\tag 1.1
$$
where $A^\perp\subset \bC T^* M$, for a subset $A\subset\bC TM$,
denotes the union over $p\in M$ of the
set of covectors at $p$ annihilating every vector in $A_p$. Real
nonvanishing sections of $T^0M$ are called {\it characteristic forms} and
sections of $T'M$ are called {\it holomorphic forms}. Thus,
characteristic forms are in particular holomorphic forms.

We shall give an alternative definition of $k_0$-nondegeneracy, as
defined in the introduction, in terms of the intrinsic geometry of
$M$. This definition appeared in [E2].
For a holomorphic form $\omega$, the Lie derivative with respect
to a CR vector field $X$ is given by $$ \Cal L_X\omega=X\lrcorner
d\omega,\tag 1.2 $$ where $\lrcorner$ denotes the interior
product, or contraction, and $d$ denotes exterior differentiation.
For $p\in M$, define the subspaces $$ T^0_pM:= E_0(p)\subset
E_1(p)\subset\ldots\subset E_k(p)\subset\ldots\subset T'_pM\tag
1.3 $$ by letting $E_k(p)$ be the linear span (over $\bC$) of the
holomorphic covectors $$ (\Cal L_{X_k}\ldots\Cal
L_{X_1}\theta)(p),\tag 1.4 $$ where $X_1,\ldots, X_k$ range over
all CR vector fields and $\theta$ over all characteristic forms
near $p$. $M$ is called {\it finitely nondegenerate} at $p\in M$
if $E_k(p)=T'_pM$ for some $k$. More precisely, we say that $M$ is
$k_0$-nondegenerate at $p$ if $$ E_{k_0-1}(p)\subsetneqq
E_{k_0}(p)=T'_pM.\tag 1.5 $$ For an argument showing that this
definition coincides with that given for embedded hypersurfaces in
the introduction, the reader is referred to [BER3] (see also [E2]). For each
$k$, set $$ F_k(p)=\bar \Cal V_p\cap E_k(p)^\perp.\tag 1.6 $$ It
was shown in [E3] that the mapping $$
(X_1,\ldots,X_k,Y,\theta)\mapsto \left<(\Cal L_{X_k}\ldots\Cal
L_{X_1}\theta)(p),Y(p)\right>,\tag 1.7 $$ defines a multi-linear
mapping $$ \Cal V_p\times\ldots\Cal V_p\times F_{k-1}(p)\times
T^0_pM \to \bC.\tag 1.8 $$ which is symmetric in the first $k$
positions. The tensor so defined for $k=1$ coincides with the
classical Levi form, and the space $F_1(0)$ is the Levi nullspace.

Let us fix a distinguished point on $M$ denoted by $0\in M$. We
choose a basis $L_{1},\ldots, L_{n}$ of the sections
$C^\infty(U,\bar \Cal V)$, where $U\subset M$ is some sufficiently
small neighborhood of $0$, adapted to the filtration $$ \bar{\Cal
V}_0=F_0(0)\supset F_1(0)\supset\ldots\supset
F_k(0)\supset\ldots\supset\{0\}\tag 1.9 $$ in the following way.
Observe that the sequence of subspaces $F_k(0)$ stabilizes at a
smallest subspace $F_{k_0}(0)$, which equals $\{0\}$ if and
only if $M$ is $k_0$-nondegenerate at $0$. Let $r_k=n-\dim_\bC
F_k(0)$ and choose $L_{1},\ldots , L_n$ so that
$L_{r_k+1}(0),\ldots L_{n}(0)$ spans $F_k(0)$ for $k=0, 1,\ldots,
k_0$. We shall use the following conventions for indices. For
$j=1,2, \ldots$, Greek indices $\alpha^{(j)},\beta^{(j)}$, etc.,
will run over the set $\{1,\ldots, r_{j-1}\}$ and small Roman
indices $a^{(j)}, b^{(j)}$, etc.,  over $\{r_{j-1}+1,\ldots, n\}$.
Capital Roman indices $A,B$, etc., will run over $\{1,\ldots,
n\}$.

Now, choose also a characteristic form $\theta$ on $M$ near
$0$. We write
$$
h_{\bar A_1\ldots\bar A_k B}:=\left<\Cal L_{\bar A_k}\ldots\Cal
L_{\bar A_1}\theta,L_{B}\right>,\tag 1.10
$$
where $\Cal L_{\bar A}:=\Cal L_{L_{\bar A}}$ and $L_{\bar
A}:=\overline{L_A}$. Note that $(h_{\bar A_1\ldots\bar A_k a^{(k)}}(0))$
represents the tensor defined by \thetag{1.7} relative to the bases
$L_{\bar A}(0)$, $L_{a^{(k)}}(0)$, and $\theta(0)$ of $\Cal V_0$, $F_k(0)$,
and $T^0_0M$, respectively.

Let $T$ be a vector field near $0$ such
that $T,L_A,L_{\bar A}$ form a basis for $C^\infty(U,\bC TM)$. Let
$\theta,\theta^A,\theta^{\bar A}$ be the dual basis for $C^\infty(U,\bC T^*M)$.
Note that, for each $k=1,\ldots, k_0$, the covectors
$\theta(0),\theta^{\alpha^{(k)}}(0)$
form a basis for $E_k(0)$. For brevity, we introduce the functions
$$
h_{\bar A_1\ldots\bar A_k}:=\left<\Cal L_{\bar A_k}\ldots\Cal
L_{\bar A_1}\theta,T\right>,\tag 1.11
$$
and also
$$
\align
 R^C_{\bar A B}:=\left<d\theta^C,L_{\bar A}\wedge
L_B\right>,&\quad R^C_{D B}:=\left<d\theta^C,L_{D}\wedge
L_B\right>\\ R^C_{\bar A}:=\left<d\theta^C,L_{\bar A}\wedge
T\right>,&\quad R^C_{B}:=\left<d\theta^C,T\wedge
L_B\right>.\tag 1.12\endalign
$$
The following identity is useful.
\proclaim{Lemma 1.13} For any nonnegative integer $k$, and
indices $A_1,\ldots, A_k,C$, $D\in \{1,\ldots, n\}$, the
following identity holds
$$
h_{\bar A_1\ldots \bar A_k\bar C D}=L_{\bar C}h_{\bar A_1\ldots\bar
A_k D}+h_{\bar A_1\ldots\bar A_k B}R^B_{\bar C D}+h_{\bar
A_1\ldots\bar A_k}h_{\bar C D}.\tag 1.14
$$
\endproclaim

\demo{Proof} Recall that $\Cal L_{\bar A_k}\ldots\Cal L_{\bar
A_1}\theta$ is a holomorphic 1-form and, by the definitions
\thetag{1.10--11}, $$ \Cal L_{\bar A_k}\ldots\Cal L_{\bar
A_1}\theta=h_{\bar A_1\ldots\bar A_k D}\theta^D+h_{\bar
A_1\ldots\bar A_k}\theta.\tag 1.15 $$ Here, and for the remainder
of this paper, we use the summation convention which states that
an index appearing in both a sub- and superscript is summed over;
e.g. $h_D\theta^D=\sum_Dh_D\theta^D$. We also have, by the
definition of the interior product, $$ h_{\bar A_1\ldots\bar
A_k\bar C D}=\left<d\Cal L_{\bar A_k}\ldots\Cal L_{\bar
A_1}\theta,L_{\bar C}\wedge L_D\right>.\tag 1.16 $$ The identity
\thetag{1.14} follows by applying the exterior derivative $d$ to
\thetag{1.15} and substituting in \thetag {1.16}. \qed\enddemo

Define $\ell_0$ to be the smallest integer $\ell$ for which
$$\left\{\aligned &h_{\bar A_1\ldots\bar A_r D}(0)=0,\quad\forall
A_1,\ldots A_r,D\in \{1,\ldots, n\},\ r<\ell\\&h_{\bar
A^0_1\ldots\bar A^0_\ell D^0}(0)\neq 0,\quad\text{for some
$A^0_1,\ldots A^0_r,D^0\in \{1,\ldots, n\}$}.
\endaligned\right.\tag 1.17$$ If no such $\ell$ exists then we set
$\ell_0=\infty$. Observe that if $M$ is $k$-nondegenerate at $0$
for some $k$, then $\ell_0\leq k$, but $\ell_0<\infty$ does not
imply finite nondegeneracy. Also, note that, for any
$r\leq\ell_0$, the subspace $F_{r-1}(0)=\bar \Cal V_0$ and, hence,
the indices $a^{r}$, $b^{r}$, etc., introduced above run over the
whole index set $\{1,\ldots,n\}$.

(Also note, by the fact that $L_A$ is adapted to the filtration
\thetag{1.9}, that if $\ell_0<\infty$ then we can take $D^0=1$ in
\thetag{1.17}.) \proclaim{Lemma 1.18} For any integer $r\geq 2$
and any integer $j\geq 0$ such that $j+r\leq \ell_0$ and indices
$A_1,\ldots,A_r,C_1\ldots, C_j,D\in\{1,\ldots,n\}$, the following
holds $$\aligned h_{\bar A_1\ldots\bar A_{r-1}\bar A_r D}(0)=&\
\big (L_{\bar A_r} h_{\bar A_1\ldots\bar A_{r-1} D}\big
)(0)\\\vdots&\\\big (L_{\bar C_1}\ldots L_{\bar C_j} h_{\bar
A_1\ldots\bar A_{r-1}\bar A_r D}\big )(0)=&\ \big (L_{\bar
C_1}\ldots L_{\bar C_j}L_{\bar A_r} h_{\bar A_1\ldots\bar A_{r-1}
D}\big )(0).\endaligned\tag 1.19$$ In particular, $$h_{\bar
A_1\bar A_2\ldots\bar A_{\ell_0} D}(0)= \big(L_{\bar
A_{\ell_0}}\ldots L_{\bar A_2} h_{\bar A_1 D}\big)(0).\tag 1.20$$
\endproclaim

\demo{Proof} The first identity in \thetag{1.19} follows
immediately by evaluating \thetag{1.14} at $0$ and using the
definition of $\ell_0$. In particular, it follows that $$\big
(L_{\bar A_r} h_{\bar A_1\ldots\bar A_{r-1} D}\big )(0)=0\tag
1.21$$ for any $2\leq r\leq \ell_0$. Now, the second identity in
\thetag{1.19} follows by applying $L_{\bar C_1}$ to
\thetag{1.14} and using \thetag{1.21}. The conclusion of Lemma
1.18 follows by induction. \qed\enddemo

Recall that $M$ is said to be of finite type at $0\in M$ if
$L_A,L_{\bar A}$ and all their repeated commutators
$$[X_m,[X_{m-1},\ldots[X_2,X_1]\ldots]],\quad  X_1,\ldots,
X_m\in\{L_1,\ldots, L_n,L_{\bar 1},\ldots, L_{\bar n}\},\tag
1.22$$ evaluated at $0$ span $\bC T_0 M$. The commutator in
\thetag{1.22} is said to have length $m$. (A commutator of length
one is simply one of the vector fields $L_A$, $L_{\bar A}$.)  If
$M$ is of finite type at $0$, then it is said to be of type $m_0$
if $m_0$ is the smallest integer for which all commutators of the
form \thetag{1.22} of lengths $\leq m_0$ span $\bC T_0 M$. Define
$\ell_1$ to be the smallest integer $\ell$ for which
$$\left\{\aligned &\left<\theta,[L_{\bar A_r},\ldots[L_{\bar
A_1},L_D]\ldots ]\right>(0) =0,\quad\forall A_1,\ldots A_r,D\in
\{1,\ldots, n\},\ r<\ell\\&\left<\theta,[L_{\bar
A^0_\ell},\ldots[L_{\bar A^0_1},L_{D^0}]\ldots ]\right>(0)\neq
0,\quad\text{for some $A^0_1,\ldots A^0_r,D^0\in \{1,\ldots,
n\}$}.\endaligned\right.\tag 1.23$$ If no such $\ell$ exists then
we set $\ell_1=\infty$. Observe that $\ell_1<\infty$ implies that
$M$ is of finite type $m_0\leq \ell_1+1$ at $0$, but the converse
is not true, i.e.\ $M$ can be of finite type at $0$ while
$\ell_1=\infty$ .

\proclaim{Proposition 1.24} If either of the two integers
$\ell_0,\ell_1$ is finite, then they are equal. Indeed, for any
$r\leq \ell_0$, it holds that $$\left<\theta,[L_{\bar
A_r},\ldots[L_{\bar A_1},L_D]\ldots ]\right>(0)=-h_{\bar
A_1\ldots\bar A_r D}(0),\tag 1.25$$ for all $ A_1,\ldots A_r,D\in
\{1,\ldots, n\}$. In particular, if $M$ is $k$-nondegenerate at
$0$, then it is also of finite type $\leq k+1$.
\endproclaim

\demo{Proof} Note that the first part of Proposition 1.24 clearly
follows from \thetag{1.25}. Hence, we shall only prove
\thetag{1.25}. For any $1$-form $\xi$ and vector fields $X$, $Y$,
we have the following well known identity (see e.g. [He])
$$\left<d\xi,X\wedge
Y\right>=X\left<\xi,Y\right>-Y\left<\xi,X\right>-
\left<\xi,[X,Y]\right>.\tag 1.26$$ Thus, for a holomorphic
$1$-form $\omega$ on $M$, we obtain $$\left<\omega,[L_{\bar A},
L_d]\right>=L_{\bar A}\left<\omega,L_D\right>-\left<\Cal L_{\bar
A}\omega,L_D\right>.\tag 1.27 $$ By applying \thetag{1.27} with
$\omega=\theta$, we deduce that $$\left<\theta,[L_{\bar
A_1},L_D]\right>=-h_{\bar A_1 D}.$$ By Lemma 1.18 and the
symmetry of the tensors $h_{\bar A_1\ldots\bar A_r D}(0)$, we then
deduce that $$\big (L_{\bar C_1}\ldots L_{\bar
C_s}\left<\theta,[L_{\bar A_1},L_D]\right>\big)(0)=-h_{\bar
C_1\ldots\bar C_s\bar A_1 D}(0),\quad \forall\, 0\leq s\leq
\ell_0-1,\tag 1.28$$ where $s=0$ in \thetag{1.28} means
$\left<\theta,[L_{\bar A_1},L_D]\right>\big(0)=-h_{\bar A_1
D}(0)$. By applying \thetag{1.27} with $\omega=\Cal L_{\bar
B_j}\ldots\Cal L_{\bar B_1}\theta$, we obtain $$\left<\Cal L_{\bar
B_j}\ldots\Cal L_{\bar B_1} \theta,[L_{\bar A_1},L_D]\right>=
L_{\bar A_1}h_{\bar B_1\ldots\bar B_j D}-h_{\bar A_1\bar B_1
\ldots\bar B_j D}.\tag 1.29$$ Hence, it follows from Lemma 1.18
and the symmetry of the $h_{\bar A_1\ldots\bar A_r D}(0)$ that
$$\left (L_{\bar C_1}\ldots L_{\bar C_s}\left<\Cal L_{\bar
B_j}\ldots\Cal L_{\bar B_1} \theta,[L_{\bar
A_1},L_D]\right>\right)(0)= 0,\quad \forall\ 1\leq j+s\leq
\ell_0-1.\tag 1.30$$ Now, assume that $$\multline\big (L_{\bar
C_1}\ldots L_{\bar C_s}\left<\theta,[L_{\bar A_r},\ldots [L_{\bar
A_1}, L_D]\right>\big)(0)=\\-h_{\bar C_1\ldots\bar C_s\bar
A_1\ldots \bar A_r D}(0),\quad \forall\, 1\leq s+r\leq
\ell_0,\endmultline\tag 1.31$$ where $s\geq 0$ and the meaning
for $s=0$ is analogous to  \thetag{1.28}, and
$$\multline\left(L_{\bar C_1}\ldots L_{\bar C_s}\left<\Cal L_{\bar
B_j}\ldots\Cal L_{\bar B_1} \theta,[L_{\bar A_r},\ldots[L_{\bar
A_1},L_D]\ldots]\right>\right)(0)= 0,\\\forall\, 1\leq j+s+r\leq
\ell_0,\endmultline\tag 1.32$$ where $j,s\geq 0$, for $r=1,\ldots
R$. Observe that we have proved this for $R=1$. Now, if
$R<\ell_0$, then the \thetag{1.31} and \thetag{1.32} follows for
all $r=1,\ldots, R+1$ by applying \thetag{1.27} and Lemma 1.18.
The verification of this is straightforward and left to the
reader. By induction, we deduce that \thetag{1.31} and
\thetag{1.32} hold for $r=1,\ldots, \ell_0$. In particular,
\thetag {1.25} holds for any $r=1,\ldots ,\ell_0$. This completes
the proof of Proposition 1.24.\qed
\enddemo

So far, everything has been done with an arbitrary choice of basis
$T,L_A,L_{\bar A}$,
except that we chose the $L_A$ to be adapted to the filtration in
\thetag{1.9} as explained above. We shall now use the fact that $M$
is embedded in $\bC^N$ and choose a particular basis.

\proclaim{Lemma 1.33} Let $M\subset\bC^N$ be a smooth real
hypersurface. Then, there is a basis $T,L_A,L_{\bar A}$ such that $T$ is
real,
the $L_A$ adapted to the filtration \thetag{1.9} as explained above, and
$$
R^C_{\bar A B}\equiv R^C_{DB}\equiv R^C_{\bar A}\equiv R^C_B\equiv 0,\tag 1.34
$$
for all indices $A,B,C,D\in\{1,\ldots,n\}$.
\endproclaim

\remark{Remark} Clearly, the conditions \thetag{1.34} are equivalent
to $d\theta^C=0$. Based on this observation, an alternative proof of
Lemma 1.33 can be given by pulling back suitable coordinate functions
from the ambient space.
\endremark

\demo{Proof} By making use of the identity \thetag{1.27}, we conclude that
$$
\aligned
R^C_{\bar A B}=-\left<\theta^C,[L_{\bar A},L_B]\right>,&\quad
R^C_{A B}=-\left<\theta^C,[L_{A},L_B]\right>,\\
R^C_{\bar A}=-\left<\theta^C,[L_{\bar A},T]\right>,&\quad
R^C_{ A}=-\left<\theta^C,[T,L_{A}]\right>.
\endaligned
\tag 1.35
$$
Hence, to prove the lemma, it is equivalent to showing that there is a
basis $T,L_A,L_{\bar A}$ with $T$ real and $L_A$ adapted to the filtration
\thetag{1.9} such that the $L_A$ commute, and $[L_{\bar A},L_B]$ and
$[L_A,T]$ are multiples of $T$. The existence of such a basis,
disregarding the adaption of the $L_A$ to the filtration, is
well known (see e.g. [BER3, Proposition 1.6.9]). Since the
adaption of the $L_A$ is a condition only at the point $0$, we may achieve
this by a linear transformation with constant coefficients to any basis $L_A$.
Such a transformation does not affect any commutator relations and, hence, the
lemma follows.\qed
\enddemo

In what follows, we shall assume that \thetag{1.34} holds.

\heading 2. A reflection identity for CR diffeomorphisms\endheading

Let $M$ be a smooth CR manifold as in the preceeding section, and
let $\hat M$ be another smooth CR manifold of the same dimension
and CR dimension, with distinguished point $\hat 0\in \hat M$. We
shall denote corresponding objects on $\hat M$ by using $\hat{}\
$; e.g.\ $\hat\Cal V\subset \bC T\hat M$ denotes the CR bundle on
$\hat M$, $\hat T,\hat L_A,\hat L_{\bar A}$ is a basis for
$C^\infty(\hat U,\bC T\hat M)$, where $\hat U$ is some
sufficiently small neighborhood of $\hat 0\in\hat M$. We shall assume
that both $M$ and $\hat M$ are embeddable, locally near $0\in M$ and
$\hat 0\in \hat M$, as real hypersurfaces in $\bC^N$. Hence,
\thetag{1.35} holds on $M$ and analogous identities on $\hat M$.

Assume that
$f\:M\to \hat M$ is a smooth CR diffeomorphism defined near $0$ in
$M$ such that $f(0)=\hat 0$. Recall that a smooth mapping $f\:M\to
\hat M$ is called CR if $f_*(\Cal V_p)\subset \hat\Cal V_{f(p)}$,
where $f_*\:\bC TM\to \bC T\hat M$ denotes the tangent mapping or push
forward, for every $p\in M$; a CR diffeomorphism is a diffeomorphism
which is CR and whose inverse is also CR. In particular, if $f$ is a
CR diffeomorphism then, for every $p\in M$ near $0$,
$f_*(\Cal V_p)= \hat\Cal V_{f(p)}$. We introduce the
smooth $GL(\bC^n)$-valued function  $(\gamma^A_B)$, and
real-valued functions $\xi,\eta^A$ so that $$
f_*(L_B)=\gamma^A_B\hat L_A,\quad f_*(L_{\bar
B})=\overline{\gamma^A_B}\hat L_{\bar A},\quad f_*(T)=\xi\hat
T+\eta^A\hat L_A+\overline{\eta^A}\hat L_{\bar A}.\tag 2.1 $$ We
can write \thetag{2.1} using matrix notation as $$
f_*(T,L_B,L_{\bar B})=(\hat T,\hat L_A,\hat L_{\bar A})\pmatrix
\xi&0&0\\\eta^A&\gamma^A_B&0\\\overline{\eta^A}&0
&\overline{\gamma^A_B}\endpmatrix .\tag 2.2 $$ By duality, we
then have $$ f^*\pmatrix \hat\theta\\\hat
\theta^A\\\hat\theta^{\bar A}\endpmatrix =\pmatrix
\xi&0&0\\\eta^A&\gamma^A_B&0\\\overline{\eta^A}&0
&\overline{\gamma^A_B}\endpmatrix \pmatrix \theta\\ \theta^B\\
\theta^{\bar B}\endpmatrix .\tag 2.3 $$

The main technical result in this section is the following, which can
be viewed as reflection identities for $\gamma^D_E$ and $\eta^D$.

\proclaim{Theorem 2.4} If $\hat M$ is $k_0$-nondegenerate at
$\hat 0\in \hat M$, then the following identities
holds for any indices $D,E\in
\{1,\ldots, n\}$,
$$
\align
\gamma^D_E= &\, r_{E}^D\big(\overline{L^J\gamma^C_A},
\overline{L^{ I}\xi};f\big),\tag 2.5
\\
\eta^D=&\, s^D\big(\overline{L^J\gamma^C_A},
\overline{L^{ I}\xi};f\big)\tag 2.6
\endalign
$$
where
$$
\aligned
&r_{E}^D\big(
\overline{L^J\gamma^C_A},
\overline{L^{ I}\xi};q\big)(p),\quad
s^D\big(\overline{L^J\gamma^C_A},
\overline{L^{ I}\xi};q\big)(p)
\endaligned\tag 2.7
$$
are smooth functions which are rational in $\overline{L^J\gamma^C_A}$ and 
polynomial in $\overline{L^{ I}\xi}$,
the indices $A$, $C$ run
over the set $\{1,\ldots, n\}$, and $J$, $I$ over all multi-indices with
$|J|\leq
k_0-1$ and $|I|\leq k_0$; here, $(p,q)\in M\times \hat M$.
Moreover, the functions in \thetag{2.7} depend only on $M$
and $\hat M$ (and not on the
mapping $f$).
\endproclaim

For the proof of Theorem 2.4, we shall make use of the following identity
$$
\left<df^*\hat \omega,X\wedge Y\right>=
\left<d\hat \omega,f_*X\wedge f_*Y\right>,\tag 2.8
$$
which holds for any 1-form $\hat \omega$ on $\hat M$ and vector fields
$X$, $Y$ on $M$. First, letting $\hat\omega=\hat\theta$, $X=L_{\bar A}$,
and $Y=L_B$, we obtain
$$
\xi h_{\bar A B}=\gamma^D_B\overline{\gamma^C_A}\hat h_{\bar C D}.\tag
2.9
$$
Here, and in what follows, we abuse the notation in the following
way. For a function $\hat c$ defined on $\hat M$,
we use the notation $\hat c$ to denote both the function $\hat c\circ
f$ on $M$ and the function $\hat c$ on $\hat M$. It should be clear from
the context which of the two functions is meant. For instance, in
\thetag{2.9}, we must have $\hat h_{\bar C D}=\hat h_{\bar C
D}\circ f$. By letting $\hat\omega=\hat\theta^E$, $X=L_{\bar A}$,
and $Y=L_B$ in \thetag{3.1}, we obtain
$$
L_{\bar A}\gamma^E_B+\eta^E h_{\bar A
B}=0.\tag 2.10
$$
Applying \thetag{2.8} with $X=L_{\bar A}$, $Y=T$,
and $\hat \omega=\hat\theta$, we obtain
$$
L_{\bar A}\xi+\xi h_{\bar A}=\xi\overline{\gamma^C_A}\hat h_{\bar
C}+\overline{\gamma^C_A} \eta^D \hat h_{\bar C D},\tag 2.11
$$
and with $\hat\omega=\hat\theta^{C}$, we obtain
$$
L_{\bar A}\eta^C+\eta^Ch _{\bar A}=0.\tag 2.12
$$
To obtain \thetag{2.11} and \thetag{2.12}, we have used the fact
$$
\aligned
\left<d\hat \omega,\hat L_{\bar C}\wedge\hat L_{\bar D}\right>= &\,\hat
L_{\bar C}\left<\hat\omega,\hat L_{\bar D}\right>-L_{\bar
D}\left<\hat\omega,\hat L_{\bar C}\right>-\left<\hat\omega,[\hat
L_{\bar C},\hat L_{\bar D}]\right> \\ =&\, 0,\endaligned
$$
which holds for any holomorphic 1-form $\hat \omega$ on $\hat M$ by
the formal integrability of the CR bundle $\hat\Cal V$. We apply
\thetag{2.8} one last time, with $\hat\omega=\hat\theta^{C}$,
$X=T$, and $Y=L_A$, to obtain
$$
T\gamma^C_A-L_A\eta^C-\eta^C \overline{ h_{\bar
A}}=0.\tag 2.13
$$

\proclaim{Lemma 2.14} For any nonnegative integer $k$, and
indices $A_1,\ldots, A_k,B,C\in \{1,\ldots, n\}$, the
following identities hold
$$
\multline
L_{\bar C}\big(\gamma^D_B\hat h_{\bar A_1\ldots\bar A_k
D}\big)=\gamma^H_B\overline{\gamma^I_C} \hat h_{\bar A_1\ldots \bar
A_k\bar I H}-\\\gamma^H_B\overline{\gamma^I_C} \hat h_{\bar A_1\ldots \bar
A_k}\hat h_{\bar I H}-\eta^H\hat h_{\bar A_1\ldots \bar
A_k H}h_{\bar C B}\endmultline\tag 2.15$$
and
$$
\multline
L_{\bar C}\big(\eta^D\hat h_{\bar A_1\ldots\bar A_k D}\big
)=\eta^H\overline{\gamma^I_C} \hat h_{\bar A_1\ldots\bar A_k
\bar I H} -\eta^H\overline{\gamma^I_C} \hat h_{\bar
A_1\ldots\bar A_k}\hat h_{\bar I H} \\-\eta^H\hat h_{\bar A_1\ldots
\bar A_k H}h_{\bar C}.
\endmultline\tag 2.16
$$
\endproclaim
\demo{Proof} We shall prove \thetag{2.15}. Recall that $\hat h_{\bar
A_1\ldots\bar A_k D}$ in
\thetag{2.15} denotes $\hat
h_{\bar A_1\ldots\bar A_k D}\circ f$ by the convention introduced in
\S 1. Hence,
$$
L_{\bar C}(\hat h_{\bar A_1\ldots\bar A_k
D})=\overline{\gamma^I_C}(\hat L_{\bar I}\hat
h_{\bar A_1\ldots\bar A_k
D}),\tag 2.17
$$
where according to our convention $\hat L_{\bar I}\hat h_{\bar
A_1\ldots\bar A_k D}=(\hat L_{\bar I}\hat h_{\bar
A_1\ldots\bar A_k D})\circ f$,
and we obtain
$$
L_{\bar C}\big(\gamma^D_B\hat h_{\bar A_1\ldots\bar A_k
D}\big)=(L_{\bar C}\gamma^D_B)\hat h_{\bar A_1\ldots\bar A_k D}+
\gamma^D_B\overline{\gamma^I_C}(\hat L_{\bar I}\hat h_{\bar
A_1\ldots\bar A_k D}).\tag 2.18
$$
Let us rewrite \thetag{2.10} as
$$
L_{\bar A}\gamma^E_B=-\eta^E h_{\bar A
B}.\tag 2.19
$$
The identity
\thetag{2.15} follows by substituting \thetag{2.19} in \thetag{2.18}
and then applying
Lemma 1.13.

The proof of the identity \thetag{2.16} is
completely analogous. Expand the left hand side by the chain rule, and
then
substitute for the derivatives of $\eta^D$ by using
\thetag{2.12}, and for the derivatives of $\hat h_{\bar A_1\ldots \bar
A_k D}$ by using Lemma 1.13.  The details are left to the reader. \qed
\enddemo

The following two lemmas will be important in establishing
Theorem 2.4.

\proclaim{Lemma 2.20} For any integer $k\geq 0$, and
indices $A_1,\ldots, A_k,B\in \{1,\ldots, n\}$, the
following identity holds
$$
\multline
\gamma^D_B\overline{\gamma^{C_1}_{A_1}}\ldots\overline{\gamma^{C_k}_{A_k}}
\hat
h_{\bar C_1\ldots\bar C_k D}=r_{\bar A_1\ldots\bar
A_k B}\big(\overline{L^J\gamma^C_A};f)+\xi \, s_{\bar A_1\ldots\bar
A_k B}\big(\overline{L^J\gamma^C_A};f)\\
+\sum_{l=1}^{k-1} \gamma^D_E\hat h_{\bar C_1\ldots\bar C_l D}
\, t^{\bar C_1\ldots\bar C_l E}_{\bar A_1\ldots\bar
A_k B}\big(\overline{L^J\gamma^C_A};f) +\sum_{l=1}^{k-1}\eta^D\hat
h_{\bar C_1\ldots\bar C_l D}
\, u^{\bar C_1\ldots\bar C_l}_{\bar A_1\ldots\bar
A_k B}\big( \overline{L^J\gamma^C_A};f),
\endmultline\tag 2.21
$$
where
$$
\aligned
& r_{\bar A_1\ldots\bar
A_k B}\big(\overline{L^J\gamma^C_A};q)(p),\quad s_{\bar A_1\ldots\bar
A_k B}\big(\overline{L^J\gamma^C_A};q)(p),\\ &t^{\bar C_1\ldots\bar
C_l E}_{\bar A_1\ldots\bar
A_k B}\big(\overline{L^J\gamma^C_A};q)(p),\quad u^{\bar C_1\ldots\bar
C_l E}_{\bar A_1\ldots\bar
A_k B}\big(\overline{L^J\gamma^C_A};q)(p)\endaligned\tag 2.22
$$
are polynomials in
$\overline{L^J\gamma^C_A}$, where $A$, $C$ run
over the indices $\{1,\ldots, n\}$ and $J=(J_1,\ldots,
J_t)\in \{1,\ldots n\}^t$ for $t\leq
k-1$, whose coefficients are smooth functions of $(p,q)\in M\times
\hat M$; here, we have used the notation  $L^J=L_{J_1}\ldots
L_{J_t}$. Moreover, the functions in \thetag{2.22} depend only on $M$
and $\hat M$ (and not on the
mapping $f$).
\endproclaim
\demo{Proof} We observe
that \thetag{2.9} satisfies the conclusion of Lemma 2.20 for
$k=1$. Assume that the
conclusion of
Lemma 2.20 holds for all integers $k=1,\ldots
j-1$. Fix indices $A_1,\ldots, A_{j-1}, B\in\{1,\ldots, n\}$, choose an
index $A_{j}\in\{1,\ldots,n\}$, and apply
$L_{\bar A_{j}}$ to \thetag{2.21} with $k={j-1}$. The statement of the
proposition for $k=j$ now follows by applying Lemma 2.14 and
substituting for $L_{\bar A_j}\xi$ using \thetag{2.11}. The proof is
completed by induction on $k$. \qed
\enddemo

\remark{Remark} In what follows, we shall use the notation $r$, $s$,
$t$, and $u$ with varying sets of sub- and superscripts for ``generic''
functions which may be different from time to time.
\endremark
\proclaim{Lemma 2.23} For any integer $k\geq 0$, and
indices $A_1,\ldots, A_k\in \{1,\ldots, n\}$, the
following identity holds
$$
\multline
\eta^D\overline{\gamma^{C_1}_{A_1}}\ldots\overline{\gamma^{C_k}_{A_k}}\hat
h_{\bar C_1\ldots\bar C_k D}=
r_{\bar A_1\ldots\bar
A_k}\big(\overline{L^J\gamma^C_A},L^{\bar I}\xi;f)\\+\sum_{l=1}^{k-1}
\gamma^D_E\hat h_{\bar C_1\ldots\bar C_l D}
\, t^{\bar C_1\ldots\bar C_l E}_{\bar A_1\ldots\bar
A_k}\big(\overline{L^J\gamma^C_A};f)
+\sum_{l=1}^{k-1}\eta^D\hat
h_{\bar C_1\ldots\bar C_l D}
\, u^{\bar C_1\ldots\bar C_l}_{\bar A_1\ldots\bar
A_k}\big( \overline{L^J\gamma^C_A};f),
\endmultline\tag 2.24
$$
where
$$
\aligned
& r_{\bar A_1\ldots\bar
A_k}\big(\overline{L^J\gamma^C_A},L^{\bar I}\xi;q)(p),\quad u^{\bar C_1
\ldots\bar
C_l E}_{\bar A_1\ldots\bar
A_k}\big(\overline{L^J\gamma^C_A};q)(p)\endaligned\tag 2.25
$$
are polynomials in
$\overline{L^J\gamma^C_A}$ and $L^{\bar I}\xi$ in the former case and in
$\overline{L^J\gamma^C_A}$ in the latter, where
$A$, $C$ run over the indices $\{1,\ldots, n\}$ and $J=(J_1,\ldots, J_t)$,
$I=(I_1,\ldots, I_{t+1})$, with $I_i,J_j\in\{1,\ldots n\}$, for
$t\leq k-1$, whose coefficients are smooth functions of $(p,q)\in
M\times \hat M$; here, we have used the notation
$L^J=L_{J_1}\ldots,L_{J_t}$ and $L^{\bar J}=L_{\bar
J_1}\ldots,L_{\bar J_t}$. Moreover, the functions in
\thetag{3.18} depend only on $M$ and $\hat M$ (and not on the
mapping $f$).
\endproclaim
\demo{Proof} We start with equation \thetag{2.11} and proceed as in
the proof of
Lemma 2.20. We leave the details to the reader.\qed
\enddemo

We are ready to prove Theorem 2.4.

\demo{Proof of Theorem $2.4$} Using the fact that the matrices $\overline
{\big(\gamma^C_A\big)}$ are invertible, we rewrite \thetag{2.21} and
\thetag{2.24} as follows 
$$
\multline
\gamma^D_B\hat
h_{\bar A_1\ldots\bar A_k D}+\sum_{l=1}^{k-1} \gamma^D_E\hat h_{\bar C_1\ldots\bar C_l D}
\, {}'t^{\bar C_1\ldots\bar C_l E}_{\bar A_1\ldots\bar
A_k B}\big(\overline{L^J\gamma^C_A};f)+\\
\sum_{l=1}^{k-1}\eta^D\hat
h_{\bar C_1\ldots\bar C_l D}
\, {}'u^{\bar C_1\ldots\bar C_l}_{\bar A_1\ldots\bar
A_k B}\big( \overline{L^J\gamma^C_A};f) 
={}'r_{\bar A_1\ldots\bar
A_k B}\big(\overline{L^J\gamma^C_A};f)+\\ {}'s_{\bar A_1\ldots\bar
A_k B}\big(\overline{L^J\gamma^C_A};f)\,\xi ,
\endmultline\tag 2.26
$$
and 
$$
\multline
\eta^D\left(\hat h_{\bar A_1\ldots\bar A_k D}+\sum_{l=1}^{k-1}\hat
h_{\bar C_1\ldots\bar C_l D}
\, {}'u^{\bar C_1\ldots\bar C_l}_{\bar A_1\ldots\bar
A_k}\big( \overline{L^J\gamma^C_A};f)\right)+\\
\sum_{l=1}^{k-1}
\gamma^D_F\hat h_{\bar C_1\ldots\bar C_l D}
\, {}'t^{\bar C_1\ldots\bar C_l F}_{\bar A_1\ldots\bar
A_k}\big(\overline{L^J\gamma^C_A};f)={}'r_{\bar A_1\ldots\bar
A_k}\big(\overline{L^J\gamma^C_A},
\overline{L^{ I}\xi};f\big).
\endmultline\tag 2.27
$$
To prove the theorem, we must show that there are $n$
choices $\underline{A}^j$, where $\underline{A}^j=A^j_1\ldots A^j_{l_j}$
with $l_j\leq k_0$, such that the linear equations
\thetag{2.26--2.27}, with $\underline{A}=\underline{A}^j$ for
$j=1,\ldots, n$ and $B=1,\ldots n$, can be solved uniquely for
$\gamma^D_B$ and $\eta^D$ near $0$. For this it suffices to show that
if, for some $(v^D_B,v^D)\in\bC^{n^2+n}$,
$$
\multline
v^D_B\hat
h_{\bar A_1\ldots\bar A_k D}(0)+\sum_{l=1}^{k-1} v^D_E\hat
h_{\bar C_1\ldots\bar C_l D} (0)
\, {}'t^{\bar C_1\ldots
\bar C_l E}_{\bar A_1\ldots\bar
A_k B}\big(\overline{L^J\gamma^C_A};f)(0)+\\
\sum_{l=1}^{k-1}v^D\hat
h_{\bar C_1\ldots\bar C_l D}(0)
\, {}'u^{\bar C_1\ldots\bar C_l}_{\bar A_1\ldots\bar
A_k B}\big( \overline{L^J\gamma^C_A};f)(0)
=0\endmultline\tag 2.28
$$
and 
$$
\multline
v^D\left(\hat h_{\bar A_1\ldots\bar A_k D}(0)+\sum_{l=1}^{k-1}\hat
h_{\bar C_1\ldots\bar C_l D}(0)
\, {}'u^{\bar C_1\ldots\bar C_l}_{\bar A_1\ldots\bar
A_k}\big( \overline{L^J\gamma^C_A};f)(0)\right)+\\
\sum_{l=1}^{k-1}
v^D_F\hat h_{\bar C_1\ldots\bar C_l D}(0)
\, {}'t^{\bar C_1\ldots\bar C_l F}_{\bar A_1\ldots\bar
A_k}\big(\overline{L^J\gamma^C_A};f)(0)=0,
\endmultline\tag 2.29
$$
for all $A_1,\ldots A_k,B\in \{1,\ldots,n\}$ and all $k\leq k_0$, then
$v^D_B=v^D=0$.  
To see this, note that \thetag{2.28--2.29}, for $k=1$, implies
directly that
$v^{\alpha^{(2)}}=v^{\alpha^{(2)}}_B=0$; recall the convention introduced in \S 1 that
the indices $\alpha^{(k+1)}$ 
run over $\{1,\ldots, r_{k}\}$, where $r_k=n-\dim F_k(0)$ as
introduced in \S 1,
and the indices $a^{(k+1})$ run over the set
$\{r_{k}+1,\ldots n\}$. Thus, since $\hat h_{A a^{(2)}}(0)=0$,
the equations \thetag{2.28--2.29} for $k=2$ reduce to  
$$
v^{a^{(2)}}_B\hat
h_{\bar A_1\bar A_2 a^{(2)}}(0)=0,\quad
v^{a^{(2)}}\hat h_{\bar A_1\bar A_2 a^{(2)}}(0)=0,\tag 2.30
$$
which in turn implies
$v^{\alpha^{(3)}}=v^{\alpha^{(3)}}_B=0$. Proceeding inductively, using
at each step 
the fact that for any
integers $1\leq j<k\leq k_0$,
$$
\hat h_{\bar A_1\ldots\bar A_j a^{(k)}}(0)=0, \tag 2.31
$$
we conclude that the equations \thetag{2.28--2.29}, for $k\leq k_0$, imply
$v^{\alpha^{(k_0+1)}}=v^{\alpha^{(k_0+1)}}_B=0$,
which is equivalent to $v^D=v^D_B=0$ since $r_{k_0+1}=n$ for a $k_0$-nondegenerate
CR manifold. This completes the proof of Theorem 2.4. \qed\enddemo

\heading 3. Proofs of Theorems 1, 2, and 3\endheading

We begin with the proof of Theorem 2. For this proof, we shall need
the following two lemmas. We shall keep the notation introduced in
previous sections.

\proclaim{Lemma 3.1} For any indices $D, E,F\in\{1,\ldots,n\}$,
multi-index
$J$, and nonnegative integer $k$, we have the following
$$
\align
L_E\overline{L^J\gamma^D_F}= &\, r^{\bar D \bar J}_{\bar F E}\big(
\overline{L^I\eta^C}\big),
\tag 3.2
\\
L_E\overline{L^JT^k\eta^D}= &s^{\bar D \bar J k}_E\big(
\overline{L^IT^m\eta^C},\overline{L^KT^{m+1}\eta^C}\big),
\tag 3.3
\endalign
$$
where the functions in \thetag{3.2--3} are smooth functions which are
rational in
the arguments appearing inside the parentheses. The indices $A$, $C$ run
over the set $\{1,\ldots, n\}$, and $I$, $K$, over all multi-indices
with $|I|\leq
|J|$, $|K|\leq |J|-1$; the integer $m$ runs from $0$ to $k$.
Moreover, the functions in \thetag{3.2--3} depend only on $M$
and $\hat M$ (and not on the
mapping $f$).\endproclaim

\demo{Proof} We shall use the following fact, which is an easy consequence
of the commutator relations established in the proof of Lemma 1.32.
For any vector field $X\in\{L_A,L_{\bar A},T\}$ and any multi-index
$J=(J_1,\ldots J_{|J|})$, we have
$$
XL^J=L^J X+\sum_{|K|
\leq |J|-1}c_{K} L^{K}T,\tag 3.4
$$
where the $c_{K}$ are smooth functions on $M$
(which depend on $X$ and $J$). To prove \thetag{3.2}, we observe that,
in view of \thetag{3.4}, we have
$$
\aligned
L_E\overline{L^J\gamma^D_F}=&\,\overline{L_{\bar E}L^J\gamma^D_F}\\=&\,
\overline{L^J L_{\bar E}\gamma^D_F}+\sum_{|K|
\leq |J|-1}c_{K} \overline{L^K T\gamma^D_F}.
\endaligned\tag 3.5
$$
The identity \thetag{3.2} follows from \thetag{2.10} and \thetag{2.13}.
The proof of \thetag{3.3} is similar, and left to the reader.\qed
\enddemo

\proclaim{Lemma 3.6} For any index $E\in\{1,\ldots,n\}$,
multi-indices
$J$ and any nonnegative integer $k$,
we have the following
$$
%\multline
L_E\overline{L^JT^k\xi}=
%\\
s_E^{\bar J k}\big(\overline{L^K\gamma^C_A},\overline{L^I T^{m}\eta^C},
\overline{L^I T^m\xi},\overline{L^K T^{m+1}\xi},T^{m}\gamma^C_A,T^m
\eta^C;f\big)
%\endmultline
\tag 3.7
$$
where the function in \thetag{3.7} is a smooth functions which are
rational in
the arguments preceding the $;$. The indices $A$, $C$ run
over the set $\{1,\ldots, n\}$, and $I$, $K$ over
all multi-indices
with $|I|\leq
|J|$, $|K|\leq |J|-1$;
the integer $m$ runs from $0$ to $k$.
Moreover, the function in \thetag{3.7} depends only on $M$
and $\hat M$ (and not on the
mapping $f$).\endproclaim

\demo{Proof} We apply
\thetag{3.4} as in the proof of Lemma 3.1 to deduce that
that $L_E\overline{L^JT^k\xi}$ is linear in $\overline{L^I T^m\xi}$,
$\overline{L^KT^{m+1}\xi}$, and $\overline{L^IT^mL_{\bar E}\xi}$. To evaluate
the latter term, we make use of \thetag{2.11} and \thetag{2.13} to deduce that
$\overline{L^JT^kL_{\bar E}\xi}$ is polynomial in $\overline{L^I T^m\xi}$,
$\overline{L^IT^m\eta^C}$, $\overline{L^K\gamma^C_A}$,
and $L^{\bar I}T^m\gamma^C_A$. Finally, we commute $L^{\bar I}$ and $T^m$
using \thetag{3.4}, and then use \thetag{2.10} to conclude that
$L^{\bar I}T^m\gamma^C_A$ is a linear function of $T^m\gamma^C_A$ and $T^m\eta^C$.
Summing up, we
obtain \thetag{3.7}. This completes the
proof of Lemma 3.6.\qed
\enddemo

The following argument is inspired by the paper [Han2]. We
shall say, for a function $u$ on $M$, that $ u\in C^a_p$
if $$ u=r(\overline{L^I\gamma^C_A},
\overline{L^IT^m\eta^C},\overline{L^IT^m\xi},L^NT^n\gamma^C_A,L^NT^n
\eta^C;f),\tag 3.8 $$ where the function in
\thetag{3.8} is a smooth functions which is rational in the
arguments preceeding the $;$. The indices $A$, $C$ run over the
set $\{1,\ldots, n\}$. The multi-indices $I$, $N$  and the
nonnegative integers $m$, $n$ run over all multi-indices with
$|I|+m\leq p$, $|N|+n\leq a$. Moreover, the function in
\thetag{3.8} should depend only on $M$ and $\hat M$ (and not on
the mapping $f$). Similarly, we shall say that $u\in
C^{a,b}_{p,q}$ if \thetag{3.8} holds with the additional
condition that $m\leq q$, $n\leq b$. Observe that by Lemma 2.30
(and the reality of $\xi$), we have
$$\gamma^D_F,\eta^C,\xi\in
C_{k_0,0}^{-1,-1},\tag 3.9
$$
where the negative ones in the superscript
signify that no terms involving $L^N\gamma^C_A$ or $L^N\eta^C$ appear. Recall
that $k_0$ is the order of nondegeneracy
of $\hat M$. By \thetag{2.10--12}, we obtain $$ L^{\bar
J}\gamma^D_F,\ L^{\bar J}\eta^C\in
C_{k_0,0}^{-1,-1}\tag 3.10 $$
for any multi-index $J$.
By
applying $L_E$ to e.g. the equations for $L^{\bar J}\gamma^D_F$
and using Lemmas 3.1 and 3.6, we conclude that $$ L_EL^{\bar
J}\gamma^D_F=r^{\bar S}_E(\overline{L^K\gamma^C_A},
\overline{L^IT^m\eta^C},\overline{L^IT^m\xi},\gamma^C_A,\eta^C;f),\tag
3.11 $$ where $|K|\leq k_0-1$, $|I|+m\leq k_0$, and $m\leq 1$. By
substituting for $\gamma^C_A$ and $\eta^C$ in \thetag{3.11} using the
equations provided by \thetag{3.10}, we conclude that $L_E\gamma^D_F\in
C^{-1,-1}_{k_0,1}$. We
obtain a similar equation for $L_EL^{\bar J}\eta^C$. Hence, we obtain
$$
L_EL^{\bar J}\gamma^D_F,\ L_EL^{\bar J}\eta^C\in
C_{k_0,1}^{-1,-1}. \tag
3.12 $$ Next, by applying
$L_{\bar F}$ to the equations for $L_EL^{\bar J}\gamma^D_F$ and
$L_EL^{\bar J}\eta^C$, provided by
\thetag{3.12}, we obtain $$ L_{\bar F}L_EL^{\bar J}\gamma^D_F,\
L_{\bar F}L_EL^{\bar J}\eta^C\in
C_{k_0+1,1}^{-1,-1}.\tag 3.13 $$ Similarly, repeated application of $L_{\bar F_1}$,
$L_{\bar F_2}$, $\ldots,L_{\bar F_k}$ yields
$$
%\multline
L_{\bar F_k}\ldots
L_{\bar F_1} L_EL^{\bar J}\gamma^D_G,\ L_{\bar F_k}\ldots L_{\bar
F_1} L_EL^{\bar J}\eta^C\in  C^{-1,-1}_{k_0+k,1}.
%\endmultline
\tag 3.14 $$
Hence, by
taking linear combinations of $L_{\bar F_k}\ldots L_{\bar
F_i}L_EL_{\bar F_{i+1}}\ldots L_{\bar F_k} L^{\bar J}$, we deduce
that
$$
[\ldots[L_E,L_{\bar F_1}],\ldots, L_{\bar F_k}]L^{\bar J}\gamma^D_G,\
[\ldots[L_E,L_{\bar F_1}],\ldots, L_{\bar F_k}]L^{\bar J}\eta^C\in C_{k_0+k,1}^{-1,-1}.
\tag 3.15 $$ Let $\ell_0\leq k_0$ be the integer
provided by Lemma 1.24 for which
$$
[\ldots[L_E,L_{\bar F_1}],\ldots, L_{\bar F_k}]=aT,\tag 3.16
$$
for some function $a$ with $a(0)\neq 0$. Then, \thetag{3.15} implies, in
particular, that $$ T\gamma^C_G,\ T\eta^C\in
C_{k_0+\ell_0,1}^{-1,-1}.\tag 3.17 $$

Before proceeding, we shall need the following result on
commutators of differential operators, which seems to be of
independent interest.

\proclaim{Proposition 3.18
} Let $\ell_0$ be
the smallest integer for which \thetag{1.23} holds. Then, for any
multi-index $J$, integer $k\geq 1$, and index $F\in\{1,\ldots, n\}$
there exist smooth
functions $a^{E_1\ldots E_m\bar F_1\ldots\bar F_s}$, $b^{E_1\ldots
E_m}_s$  such that
$$\multline \sum_{m=1}^{|J|+k}\sum_{s=0}^{m\ell_0} a^{E_1\ldots
E_m\bar F_1\ldots\bar F_s} [\ldots [L_{E_1}\ldots L_{E_m},L_{\bar
F_1}],L_{\bar F_2}]\ldots,L_{\bar F_s}]=L^{J}T^k.\endmultline\tag
3.19 $$
and
$$
%\multline
\sum_{m=1}^{|J|+k}\sum_{s=0}^{k} b_s^{E_1\ldots
E_m} \underbrace{[\ldots [L_{E_1}\ldots L_{E_m},L_{\bar
F}],L_{\bar F}]\ldots,L_{\bar F}]}_{\text{\rm length $s$
}}=(h_{\bar F
1})^{p}L^{J}T^k,
%\endmultline
\tag
3.20 $$
where $p=k+|J|-|J|_1+1$ and $|J|_1$ denotes the number of occurences
of the index $1$ in the multi-index $J$; here, the
length of the commutator $[\ldots [X,Y_{1}],Y_{2}]\ldots,Y_{s}]$ is
$s$.
\endproclaim
\demo{Proof} In this proof, we shall use the following notation to
simplify the notation, $$ C_{E_1\ldots E_m,\bar F_1\ldots\bar
F_s}:=[\ldots[[L_{E_1}\ldots L_{E_m},L_{\bar F_1}],L_{\bar
F_2}]\ldots,L_{\bar F_s}],$$ where $C_{E_1\ldots E_m}$ is
understood to mean $L_{E_1}\ldots L_{E_m}$. Using bilinearity of
the commutator and the identity $$ [AB,C]=A[B,C]+[A,C]B,\tag 3.21
$$ a straightforward induction shows that $$\multline
[\ldots[[L_{E_1}L_{E_2}\ldots L_{E_m},L_{\bar F_1}],L_{\bar
F_2}]\ldots,L_{\bar F_s}]=\\
\sum_{(\underline{i},\underline{j})\in P_2(s)}C_{E_1,\bar
F_{i_1}\ldots L_{\bar F_{i_{s-l}}}}C_{E_2\ldots E_m,\bar
F_{j_1}\ldots \bar F_{j_l}},\endmultline\tag 3.22 $$ where
$P_2(s)$ denotes the set of all partitions of $\{1,\ldots,s\}$
into two disjoint increasing sequences $\underline{i}=(i_1,\ldots,
i_{s-l})$, $1\leq i_1<\ldots i_{s-l}\leq s$, and
$\underline{j}=(j_1,\ldots,j_{l})$, $1\leq j_1<\ldots<j_l\leq s$
for $l=0,\ldots s$. (Of course, for e.g. $l=0$ the partition is
understood to be the trivial one $\underline{i}=(1,\ldots, s)$ and
$\underline{j}=\emptyset$.) Similarly, if we denote by $P_m(s)$
the set of all partitions of $\{1,\ldots,s\}$ into $m$ disjoint
increasing sequences $\underline{i^t}=(i^t_1,\ldots, i^t_{s_t})$,
$1\leq i^t_1<\ldots i^t_{s_t}\leq s$, $t=1,\ldots m$, and $\sum
s_t=s$ (allowing empty sequences), then we have $$ \multline
[\ldots[[L_{E_1}L_{E_2}\ldots L_{E_m},L_{\bar F_1}],L_{\bar
F_{2}}]\ldots,L_{\bar F_s}]=\\
\sum_{(\underline{i^1},\ldots,\underline{i^m})\in
P_m(s)}C_{E_1,\bar F_{i^1_1}\ldots L_{\bar F_{i^1_{s_1}}}}\ldots
C_{E_m,\bar F_{i^m_1}\ldots \bar
F_{i^m_{s_m}}},\endmultline\tag 3.23 $$ Observe that $C_{E,\bar
F_1\ldots \bar F_s}=a_{E,\bar F_1\ldots\bar F_s}T$, for some function
$a_{E,\bar F_1\ldots\bar F_s}$ such that $$ a_{E,\bar
F_1\ldots\bar F_s}(0)=0,\quad \forall s<\ell_0,\quad\text{\rm
and}\quad a_{E,\bar F_1\ldots\bar F_{\ell_0}} (0)=h_{\bar
F_1\ldots\bar F_{\ell_0} E}(0)\neq 0,$$ for some choice of
$F_1,\ldots, F_{\ell_0}$. Hence, with $s=\ell_0$ we obtain, by
\thetag{3.23} and Lemma 1.24, $$ \multline [\ldots
[L_{E_1}\ldots L_{E_m},L_{\bar F_1}],L_{\bar F_2}]\ldots,L_{\bar
F_{\ell_0}}]=\sum_{l=1}^m (h_{\bar F_1\ldots \bar F_{\ell_0} E_l}+o(1))
L_{E_1}\ldots \widehat {L_{E_l}} \ldots L_{E_m} T
\\
+
\sum_{{|K|+p= m}\atop{|K|\leq m-2}} b_{Kp}L^KT^p+
\sum_{|K|\leq m-2} c_{K}L^KT,
\endmultline
\tag 3.24 $$ where $b_{kp}(0)=0$, $\widehat{L_{E_l}}$ means
that factor is omitted, and $o(1)$ denotes a function vanishing at
$0$; the last sum in \thetag{3.24} arises from
arranging (by commuting) so that the vector field $T$ comes last
in the first sum. Recall, from \S 1, that for each index
$\alpha^{(\ell_0)}\in \{1,\ldots,r_{\ell_0}\}$ there exist
$F_1,\ldots,F_{\ell_0}$ so that $h_{\bar F_1\ldots\bar F_{\ell_0}
\alpha^{(\ell_0)}}(0)\neq 0$. For this argument, we only need the
fact that there exist $F_1,\ldots,F_{\ell_0}$ so that $h_{\bar
F_1\ldots\bar F_{\ell_0} 1}(0)\neq 0$. We choose
$F_1,\ldots,F_{\ell_0}$ to be minimal, in the
lexicographical ordering ($A_1\ldots A_s<B_1\ldots B_s$ if, for some
$r\leq s$, $A_i\leq B_i$ for $i<r$, and $A_r<B_r$), with this property. 
Setting $E_1=\ldots E_m=1$,
we observe that we can solve for $L_1^{m-1}T$ in \thetag{3.24}.
Setting $E_1=\ldots E_{m-1}=1$ and $L_{E_m}=L_E$ with $E\geq 2$,
we can then solve for $L_1^{m-1}L_ET$. Proceeding inductively, we
see that we can solve for any $L^JT$, with $|J|=m-1$, in terms of
$$ [\ldots [L_{E_1}\ldots L_{E_m},L_{\bar F_1}],L_{\bar
F_2}]\ldots,L_{\bar F_{\ell_0}}], \quad b_{Kp} L^KT^p, \quad L^K
T, $$ where $K$ runs over multi-indices with $|K|\leq m-2$, and
$p$ over positive integers such that $|K|+p=m$, and each
$b_{kp}(0)=0$. Next, letting $s=2\ell_0$ and $F_{\ell_0+l}=F_{l}$
for $l=1,\ldots \ell_0$, we obtain (by also using that $h_{\bar
F_1\ldots\bar F_s E}$ is symmetric in the first $s$ indices
$F_1,\ldots F_s$) 
$$ 
\multline [\ldots
[L_{E_1}\ldots L_{E_m},L_{\bar F_1}],L_{\bar F_2}]\ldots,L_{\bar
F_{2\ell_0}}]=\sum_{l=1}^m a_{E_j\bar F_1\ldots \bar F_{2\ell_0}}
L_{E_1}\ldots \widehat {L_{E_l}} \ldots L_{E_m} T
\\
+c_{\bar F}\sum_{1\leq l_1<l_2\leq m} (h_{\bar F_1\ldots \bar
F_{\ell_0} E_{l_1}}h_{\bar F_1\ldots \bar
F_{\ell_0} E_{l_2}}+o(1)) L_{E_1}\ldots
\widehat {L_{E_{l_1}}}
\ldots \widehat {L_{E_{l_2}}}\ldots L_{E_m} T^2\\+
\sum_{{|K|+p= m}\atop{|K|\leq m-3}} o(1)L^KT^p+
\sum_{{|K|+p\leq m-1}\atop {p=1, 2}} c_{Kp}L^KT^p,
\endmultline
\tag 3.25 $$ where $c_{\bar F}$ is some combinatorial factor ($>0$)
which depends on the minimal set of indices
$F_1,\ldots,F_{\ell_0}$. Using the fact that we have 
already solved for the $L^JT$, 
$|J|=m-1$, in terms of $b_{K2}L^K T^2$ where $b_{K2}(0)=0$, a
similar argument to the one used above shows that we can solve for
each $L^{J}T^2$, $|J|=m-2$, in terms of $$ [\ldots [L_{E_1}\ldots
L_{E_m},L_{\bar F_1}],L_{\bar F_2}]\ldots,L_{\bar F_{2\ell_0}}],
\quad o(1) L^KT^p, \quad L^Q T, \quad L^KT^2 $$ where $K$, $Q$,
runs over multi-indices with $|K|\leq m-3$, $|Q|\leq m-2$, and $p$
over positive integers such that $|K|+p=m$.  Proceeding by
induction over $k$ (with the total order $m$ fixed), we eventually
find that we can solve completely for $T^m$ in terms of $$ [\ldots
[L_{E_1}\ldots L_{E_m},L_{\bar F_1}],L_{\bar F_2}]\ldots,L_{\bar
F_{m\ell_0}}], \quad L^KT^p,$$ with $|K|+p\leq m-1$. Substituting
back, we obtain $$\multline \sum_{n=0}^{k\ell_0} a^{E_1\ldots
E_m\bar F_1\ldots\bar F_n} [\ldots [L_{E_1}\ldots L_{E_m},L_{\bar
F_1}],L_{\bar F_2}]\ldots,L_{\bar F_n}]=L^{J}T^k\\+\sum_{|K|+p\leq
m-1} c_{Kp}L^KT^p,\endmultline $$ where $m=|J|+k$. The proof of
\thetag{3.19} is completed by a simple induction on the total
degree $m$.

For the proof of \thetag{3.20}, we proceed analogously by first
setting $s=1$ and $E_1=\ldots=E_m=1$ in \thetag{3.23}. We find that
$$
mh_{\bar F 1}L_1^{m-1}T=C_{E_1\ldots E_m,\bar
F}+\sum_{|K|\leq m-2} c_K L^K T.\tag 3.26
$$
Next, with $E_1=\ldots =E_{m-1}=1$ and $E_m=E$, we obtain
$$
(m-1)h_{\bar F 1}L_1^{m-2}L_ET+h_{\bar F E}L_1^{m-1}T=C_{E_1\ldots E_m,\bar
F}+\sum_{|K|\leq m-2} c'_K L^K T.\tag 3.27
$$
Thus, multiplying by $h_{\bar F 1}$ and using \thetag{3.26}, we
obtain \thetag{3.20} for a multi-index $J=(1,\ldots,1,E)\in\{1,\ldots, n\}^{m-1}$
and $k=1$. Similarly, we obtain \thetag{3.20} for
arbitrary multi-indices $J$ and $k=1$. Proceeding inductively, setting
$s=2,3,\ldots k$, using \thetag{3.23}, and multiplying through by a
suitable power of $h_{\bar F 1}$ to apply the results obtained in
previous steps, we arrive at \thetag{3.20}. The details are left to
the reader. This completes the proof of
Proposition 3.18.
\qed
\enddemo

To complete the proof of Theorem 2, let us observe the following
schematic diagram which describes the action of applying the
operators $L_E$ to elements in $C^{-1,-1}_{q+k,q}$ $$ \multline
C^{-1,-1}_{q+k,q}@> L_{E_1}>> C^{q,q}_{q+k,q+1} @>L_{E_2}>>
C^{q+1,q+1}_{q+k,q+2} @>L_{E_3}>> \ldots
@>L_{E_k}>>\\@>L_{E_k}>> C^{q+k-1,q+k-1}_ {q+k,q+k}@> L_{E_{k+1}}>>
C^{q+k,q+k}_{q+k,q+k}@>L_{E_{k+2}}>>
C^{q+k+1,q+k}_{q+k,q+k}@>L_{E_{k+3}}>>\ldots
\endmultline\tag 3.28
$$ The verification of the diagram is straightforward using Lemmas
3.1 and 3.6, and the details are left to the reader. Similarly,
we have $$ C^{a+k,a}_{p,q}@> L_{\bar F_1}>> C^{a+k,a+1}_{p+1,q}
@>L_{\bar F_2}>>\ldots @>L_{\bar F_k}>> C^{a+k,a+k}_{p+k,q}@>
L_{\bar F_{k+1}}>>C^{a+k,a+k}_{p+k+1,q}@>
L_{\bar F_{k+2}}>>\ldots\tag 3.29 $$ We claim that the following
holds for any multi-index $J$ and nonnegative integer $k$,
$$L^JT^k\gamma^D_F,\ L^JT^k\eta^D\in
C^{-1,-1}_{k_0+m\ell_0,\min(k_0,m)}, \tag 3.30 $$ where
$m=|J|+k$. Observe that \thetag{3.30} holds for $m=1$ by
\thetag{3.12} and \thetag{3.17}. We shall prove \thetag{3.30}
by induction on $m$. Thus, assume that \thetag{3.30} holds for
all $m\leq m_0-1$. By Proposition 3.18, we can produce the
differential operator $L^J T^k$ by taking linear combinations of
operators of the form $L^{\bar P} L^Q L^{\bar R}$, where
$|P|+|R|\leq m_0\ell_0$, $|Q|\leq m_0$, and $m_0=|J|+k$. Applying
first $L^QL^{\bar R}$ to e.g. $\gamma^D_F$ we conclude, using
\thetag{3.10} and the diagram \thetag{3.28}, that $L^QL^{\bar
R}\gamma^D_F\in C^{m_0-1,\min(k_0,m_0-1)}_{k_0,\min(k_0,m_0)}$.
By applying $L^{\bar P}$ to the equation for $L^Q L^{\bar
R}\gamma^D_F$ and using the diagram \thetag{3.29}, we obtain
$L^JT^k\gamma^D_F\in
C^{m_0-1,m_0-1}_{k_0+m_0\ell_0,\min(k_0,m_0)}$.
The
conclusion $L^JT^k\gamma^D_F\in
C^{-1,-1}_{k_0+m_0\ell_0,\min(k_0,m_0)}$ follows by using the
induction hypothesis to substitute for the $L^I T^m\gamma^C_A$ and
$L^I T^m\eta^C$, with $|I|+m\leq m_0-1$, that appear in the
equation for $L^JT^k\gamma^D_F$.
By applying
the same argument to $\eta^D$, we conclude that \thetag{3.30}
holds for $m= m_0$ and, hence, for all $m$ by induction. This
proves the claim. In particular, we then have $$ T^k\gamma^D_F,
T^k\eta^D\in C^{-1,-1}_{k_0+k\ell_0,\min(k_0,k)}. \tag 3.31 $$ By
applying $L^J$ to these equations and using \thetag{3.28}, we
deduce that $$ L^JT^k\gamma^D_F,\ L^JT^k\eta^D\in
C^{a(J,k),q'(J,k)}_{k_0+k\ell_0,q(J,k)},\tag 3.32 $$ where $$
\align a(J,k) =&\, \min(k_0,k)+|J|-1\tag 3.33\\ q(J,k) =&\,
\min(k_0+k\ell_0,\min(k_0,k)+|J|),\tag 3.34\endalign
$$
and
$$ 
q'(J,k) =
\min(k_0+k\ell_0,\min(k_0,k)+|J|-1).
$$
Observe that \thetag{3.32} implies
$$
L^JT^k\gamma^D_F,\ L^JT^k\eta^D \in
C^{|J|+k-1,k_0+k_0\ell_0}_{k_0+k_0\ell_0,k_0+k_0\ell_0},\quad \forall k\:0\leq
k\leq k_0.\tag 3.35
$$
Now, from \thetag{3.4} it follows that ${L^RT^k\xi}=T^kL^R\xi$ modulo terms of
the form $T^k L^S\xi$ with
$|S|<|R|$. Hence, by applying \thetag{2.11--13}, we conclude that
$L^RT^k\xi$ is a polynomial in $T^m\xi=\overline{T^m\xi}$,
$\overline{T^m\eta^C}$, ${L^K\gamma^C_A}$, and
${L^IT^{m-1}\eta^C}$, where $m\leq k$, $|K|\leq |R|-1$, and
$|I|\leq |R|$. It follows that we also have
$$
L^JT^k\xi \in
C^{|J|+k-1,k_0+k_0\ell_0}_{k_0+k_0\ell_0,k_0+k_0\ell_0},\quad \forall k\:0\leq
k\leq k_0.\tag 3.36
$$
Let us introduce the new class $D^{a,b}_{p,q}$
consisting of functions $u$ for which there is an equation
$$ u=r(\overline{L^I\gamma^C_A},
\overline{L^IT^m\eta^C},L^NT^n\xi,L^NT^n\gamma^C_A,L^NT^n
\eta^C;f),\tag 3.37 $$
where $|I|+m\leq p$, $m\leq q$, $|N|+n\leq a$, $n\leq b$, and the
function $r$
is rational in the arguments preceeding the ``;'' and only depends on
$M$ and $\hat M$. By the
above remarks concerning $L^RT^k\xi$ and \thetag{3.31}, we have
$$
T^k\gamma^D_F,\ T^k\eta^D\in
D^{\min(k_0,k),\min(k_0,k)}_{k_0+k\ell_0,\min(k_0,k)}.\tag 3.38
$$
Since equations of the form \thetag{3.37} do not involve terms of the
form $\overline{L^IT^m\xi}$, we obtain a different diagram describing
the action of $L_E$ on the
classes $D^{a,b}_{p,q}$, namely
$$
\multline
D^{a,b}_{q+k,q}@> L_{E_1}>>
D^{a+1,b}_{q+k,q+1} @>L_{E_2}>>\ldots @>L_{E_k}>> D^{a+k-1,b}_{q+k,q+k}@>
L_{E_{k+1}}>>\\@>
L_{E_{k+1}}>>
D^{a+k,b}_{q+k,q+k}@>L_{E_{k+2}}>> D^{a+k+1,b}_{q+k,q+k}@>L_{E_{k+3}}>>\ldots
\endmultline\tag 3.39
$$
By \thetag{3.38--39}, we deduce
$$
L^JT^k\gamma^D_F,\ L^J T^k\eta^D\in
D^{k_0+|J|,k_0}_{k_0+k\ell_0,\min(k_0+k\ell_0,k_0+|J|)},\quad\forall
k\:k \geq k_0+1.\tag 3.40
$$
We have the following technical, but important, lemma.
\proclaim{Lemma 3.41} For any multi-index $J$, and nonnegative
integer $k\leq k_0$, we have
$$
L^JT^k\gamma^D_F,\ L^JT^k\eta^D,\ L^JT^k\xi\in
C^{-1,-1}_{k_0+(k_0+k_0\ell_0)\ell_0,
k_0+(k_0+k_0\ell_0)\ell_0}.\tag 3.42
$$
\endproclaim
\demo{Proof} By \thetag{3.35--36}, we have $$ L^JT^k\gamma^D_F,\
L^JT^k\eta^D,\ L^JT^k\xi\in
C^{|J|+k-1,k_0+k_0\ell_0}_{k_0+k_0\ell_0,k_0+k_0\ell_0}.\tag
3.43 $$ Observe that \thetag{3.43} reduces the total order of the
unconjugated terms by at least one. Now, in the equations given by
\thetag{3.43}, there may appear terms of the form
$L^{I^1}T^{k_1}\gamma^C_A$, $L^{I^1}T^{k_1}\eta^C$, where
$|I^1|+k_1\leq |J|+k-1$, and $k_1\leq k_0+k_0\ell_0$. For those
term with $k_0+1\leq k_1\leq k_0+k_0\ell_0$, we may apply
\thetag{3.40} to deduce that $$ L^{I^1}T^{k_1}\gamma^C_A,\
L^{I^1}T^{k_1}\eta^C\in
D^{k_0+|I^1|,k_0}_{k_0+(k_0+k_0\ell_0)\ell_0,\min(k_0+(k_0+k_0\ell_0)\ell_0,
k_0+|I^1|)}.\tag
3.44 $$ Note that, since $k_1\geq k_0+1$ and $|I^1|+k_1\leq
|J|+k-1$, we have $k_0+|I^1|\leq |J|+k-2$. For those terms
$L^{I^1}T^{k_1}\gamma^C_A$, $L^{I^1}T^{k_1}\eta^C$ with
$k_1\leq k_0$, we may apply \thetag{3.35} again. In any case, we have
reduced the total order of the unconjugated terms by two. In the
equations given by \thetag{3.44}, there may appear terms of the
form $L^{I^2}T^{k_2}\gamma^C_A$, $L^{I^2}T^{k_2}\eta^C$, and also
$L^{I^2}T^{k_2}\xi$, where $|I^2|+k_2\leq k_0+|I_1|\leq |J|+k-2$,
and $k_2\leq k_0$. We again substitute for these terms, using the
equations given by \thetag{3.35--36}. This reduces the total
order of the unconjugated terms another step. Proceeding in this
way, alternately substituting using either \thetag{3.35--36} or
\thetag{3.40}, we eliminate all the unconjugated terms (in a
finite number of steps). At each step we introduce new conjugated
terms, but in view of \thetag{3.35--36} and \thetag{3.40}, the
total order of these never exceed $k_0+(k_0+k_0\ell_0)\ell_0$.
This completes the proof of Lemma 3.41.\qed
\enddemo

By substituting, using Lemma 3.41, for the conjugated terms
that appear in the equations provided by \thetag{3.30},
we conclude that for any multi-index $J$
and nonnegative integer $k$, we have
$$
L^JT^k\gamma^D_F,\ L^JT^k\eta^D\in
D^{k_0+(k_0+k_0\ell_0)\ell_0,k_0+(k_0+k_0\ell_0)\ell_0}_{-1,-1}.\tag
3.45
$$
By using \thetag{2.10} and \thetag{2.12}, we
conclude that for any multi-indices $R$ and $S$, any nonnegative
integer $k$, and any indices $D,F\in \{1,\ldots n\}$,  there are
smooth functions, which are rational in their
arguments preceeding the ``;'', such that
$$
\aligned
L^RT^kL^{\bar S} \gamma^D_F =&\, r^{R\bar S k}
\big(L^IT^j\gamma^C_A,L^IT^j\eta^C,L^IT^j\xi;f\big),
\\
L^RT^kL^{\bar S} \eta^D_F =&\, s^{R\bar S k}
\big(L^IT^j\gamma^C_A,L^IT^j\eta^C,L^IT^j\xi;f\big),
\endaligned\tag 3.46
$$ where $|I|+j\leq k_0+(k_0+k_0\ell_0)\ell_0$. Finally, by
using \thetag{2.11}, its complex conjugate, and Proposition
3.18, it is not difficult to see that $L^RT^kL^{\bar S}\xi$ can
be expressed in terms of $\xi$ and derivatives of $\gamma^C_A$,
$\overline{\gamma^C_A}$, $\eta^C$, and $\overline{\eta^C}$. Thus,
in view of \thetag{3.46}, we also have, for any $R$, $S$, and
$k$, $$ L^RT^kL^{\bar S} \xi= t^{R\bar S k}
\big(L^IT^j\gamma^C_A,L^IT^j\eta^C,L^IT^j\xi,\overline{L^IT^j\gamma^C_A},
\overline{L^IT^j\eta^C}, \overline{L^IT^j\xi} ;f\big),\tag 3.47
$$ where $I$ and $j$ run over the same indices as in
\thetag{3.46}. Now, recall that $\ell_0\leq k_0$. The conclusion
of Theorem 2 follows by writing \thetag{3.46--47}, for all $R$,
$S$, $k$ such that $$|R|+|S|+k=k_0+(k_0+k_0^2)k_0+1$$ in any
coordinate systems $x=(x_1,\ldots,x_{2n+1})$ and $\hat x=(\hat
x_1,\ldots\hat x_{2n+1})$ for $M$ and $\hat M$ near the points
$0\in M$ and $\hat 0\in \hat M$, respectively, and observing that the
same system of differential equations holds for any CR mapping $f$
sending a neighborhood of $0$ in $M$ into $\hat M$ with $f(0)$
sufficiently close to $\hat 0$. This completes the
proof of Theorem 2.\qed

\remark{Remark} We would like to point out that a much simpler conclusion of the
proof of Theorem 2 can be given in the case $\ell_0=1$, i.e.\ when
the Levi form of
$M$ has at least one nonzero eigenvalue at $0$. We can then use the
commutator identity \thetag{3.20} instead of \thetag{3.19} to conclude
$$L^JT^k\gamma^D_F,\ L^JT^k\eta^D\in
C^{-1,-1}_{k_0+k,\min(k_0,m)}, \tag 3.48 $$
instead of \thetag{3.30}. By substituting for conjugated
terms, using only \thetag{3.48}, we obtain directly equations of the form
\thetag{3.46} in which $|I|+j\leq 2k_0$. We invite the reader to carry
out the details in this case. Observe that the system of differential
equations obtained for $f$ using this argument is of order $2k_0+2$
rather than $k_0^3+k_0^2+2$ as given by Theorem 2 (or
$k_0+(k_0+k_0\ell_0-1)\ell_0+2=3k_0+1$ for $\ell_0=1$, which
is the order that actually follows from the proof of Theorem 2
presented above).

A similar simpification of the proof in the general case would be
possible if one could prove that it suffices to take the sum over $s$
in \thetag{3.19} to run from $s=0$ to $s=k\ell_0$ instead of all the
way up to $s=m\ell_0$. 
\endremark

\demo{Proof of Theorem $1$} The system of differential equations
\thetag{0.3} is a so-called complete system of order $k_0^3+k_0^2+k_0+2$.
In particular, any solution is completely determined by
its $(k_0^3+k_0^2+k_0+1)$-jet at $0\in M$ (see e.g.\ [BCG$^3$]. cf. also
[Han2, Proposition 2.2]). On the other hand, if $x\to Z(x)$
is an embedding of $M$ into $\bC^N$ sending $p_0\in M$ to $0\in
\bC^N$ and $x'\to Z'(x)$ is an embedding of $M'$ sending $p_0'$ to
$0\in \bC^N$, then for any smooth CR mapping $f\:M\to \hat M$,
with $f(p_0)=p_0'$,
there exists  (see e.g.\
[BER3, Proposition 1.7.14]) a formal power series mapping $Z'=F(Z)$,
with $F(0)=0$,
sending $M$ into $M'$ (i.e. $\rho(Z,\bar Z)$ divides
$\rho'(F(Z),\overline{F(Z)})$ in the ring of formal power series in
$Z,\bar Z$; cf.\ e.g.\ [BER4]) such that
$$Z'(f(x))\sim F(Z(x)),\tag 3.49$$ where $\sim$ denotes equality as
formal power series. Also, by [BER4, Theorem 2.1.1], the
$2k_0$-jet at $0$ of any invertible formal mapping $Z'=F(Z)$, with
$F(0)=(0)$, sending $M$ into $M'$ determines the series $F(Z)$
completely. In particular, it
follows from \thetag{3.49} that the
$2k_0$-jet at $p_0$ of a CR diffeomorphism $f\:M\to M'$, with $f(p_0)=p_0'$,
determines its $(k_0^3+k_0^2+k_0+1)$-jet at $p_0$. Hence, the conclusion of Theorem
1 follows from Theorem 2. \qed  \enddemo

\demo{Proof of Theorem $3$} We shall need the following proposition.

\proclaim{Proposition 3.50} If a smooth real hypersurface
$M\subset\bC^N$ is holomorphically nondegenerate at $p_0\in M$, then
there exists an open neighborhood $U$ of $p_0\in M$ and a dense open
subset $U'\subset U$ such that $M$ is $(N-1)$-nondegenerate at every
$p\in U'$.
\endproclaim

\demo{Proof} The statement that, under the hypotheses in the
proposition, there exists an open neighborhood $U$ of
$p_0$ such that $M$ is finitely nondegenerate on a dense open subset
$U''\subset U$ is a consequence of [BER3, Theorem 11.7.5 (iii)]. To
prove Proposition 3.50, it suffices to show that if $M$ is not
$k$-nondegenerate, for any $k\leq N-1$, on an open set $V$, then
$M$ is in fact not finitely nondegenerate at any $p\in V$. Recall the
subspaces $E_j(p)\subset T'_pM$ defined for $j=0,1,\ldots$ in \S
1. Assume that $E_{N-1}(p)$ is a proper subspace of $T'_pM$ for
every $p\in V$, i.e.\ $M$ is not
$k$-nondegenerate, for any $k\leq N-1$, in $V$. Since
$\dim_\bC T'_p M=N$, we conclude,
by \thetag{1.3}, that there must be an
open subset $V'\subset V$ and an integer $1\leq \ell\leq N-1$ such
that
$$
E_{\ell-1}(q)=E_{\ell}(q),\quad\forall q\in V'.\tag 3.51
$$
We claim that if
$E_{\ell-1}(q)=E_{\ell}(q)$ for all $q$ in some open sufficiently
small set $V'\subset M$, then in
fact $E_{\ell-1}(q)=E_k(q)$ for all $k\geq \ell$ and all $q\in
V'$. To see this, observe that \thetag{3.51} implies
that, for every $A_1,\ldots,A_\ell\in\{1,\ldots, N\}$, there are
smooth functions $a^{\bar C_1\ldots\bar C_j}_{\bar
A_1\ldots \bar A_\ell}$ such that
$$
\Cal L_{\bar A_{\ell}}\ldots \Cal L_{\bar
A_1}\theta=\sum_{j=0}^{\ell-1} a^{\bar C_1\ldots\bar C_j}_{\bar
A_1\ldots \bar A_\ell} \Cal L_{\bar C_j}\ldots\Cal L_{\bar
C_1}\theta\tag 3.52
$$
in $V'$. \thetag{3.52} implies that $E_{\ell+1}(q)=
E_{\ell}(q)$ for $q\in V'$, and the claim follows by induction. We
conclude that $M$ is not finitely nondegenerate in $V'$. A simple
connectedness argument applied to each component of $V$ proves that
$M$ cannot be finitely nondegenerate at any point in $V$. This
completes the proof of Proposition 3.50.\qed
\enddemo

We return to the proof of Theorem 3. The fact that $M$ is minimal at
$p_0$ implies, by a theorem of Trepreau (see e.g. [BER3, Theorem
8.1.1]; the analogous result in higher codimensions was proved by
Tumanov), that for any open neighborhood $U$ of $p_0$ in $M$,
there exists an open connected set $\Omega\subset \bC^N$ (on ``one
side of $M$'') such that
$U':=\overline{\Omega}\cap M\subset U$ is an open neighborhood of
$p_0$ and every smooth CR function in $U$ is the smooth boundary value
of a holomorphic function in $\Omega$. We deduce by the uniqueness of
boundary values of holomorphic functions, Proposition 3.50, and
Theorem 1 that there exists $p_1\in U'$ such that if $f_1,f_2\:U\to
M'$ are smooth CR diffeomorphisms into some smooth real hypersurface
$M'\subset \bC^N$ and $\partial^\alpha f_1(p_1)=\partial
^\alpha f_2(p_1)$, for all $|\alpha|\leq 2(N-1)$, then $f_1=f_2$ in
$U'$. Using this fact, the proof of Theorem 3 is completed exactly as
the proof of [BER2, Theorem 2]. Choose $Y_1,\ldots, Y_m\in
\aut(M,p_0)$ which are linearly independent over $\bR$, and denote by
$F(x,y)$, where $x=(x_1,\ldots, x_{2N-1})$ is some local coordinate
system on $M$ near $p_0$ and $y=(y_1,\ldots,y_m)\in \bR^m$,  the
time-one map of the flow $\exp t(y_1Y_1+\ldots+y_m Y_m)$, for $y$
in some sufficiently small neighborhood $V$ of the origin $\bR^m$. The
arguments in [BER2] combined with the uniqueness result stated above,
for a suitably chosen open neighborhood of $U$ of $p_0$ in $M$, shows
that the mapping $V\to J^{2(N-1)}(M,M')_{p_1}$, given by
$$
y\mapsto (\partial_x^\alpha F(p_1,y))_{|\alpha|\leq 2(N-1)},\tag 3.53
$$
is smooth and injective. Hence, $m\leq \dim_\bR
J^{2(N-1)}(M,M')_{p_1}$ which proves Theorem 3. \qed
\enddemo

\demo{Proof of Theorem $4$} The conclusion of Theorem 4 is a direct consequence of
Theorem 2 and Proposition 3.54 below. We shall use the notation
$J^k(\bR^q,\bR^m)_0$ 
for the space of $k$-jets at $0\in \bR^q$ of smooth mappings
$f\:\bR^q\to \bR^m$, and $\lambda^k=(\lambda^\beta_i)$, $|\beta|\leq
k$ and $i=1,\ldots, m$, for the natural coordinates on
this space in which the $k$-jet of $f$ is given by
$\lambda^\beta_i=\partial^\beta_xf_i(0)$. 

\proclaim{Proposition 3.54} Let $U\subset
J^{k}(\bR^q,\bR^m)_{0}\times \bR^q$ be an open domain. Let  $r^\alpha_j(\lambda^{k})(x)$,
for any multi-index $\alpha\in\Bbb
Z_+^{m}$ with $|\alpha|=k+1$ and any
$j=1,\ldots,m$, be smooth ($C^\infty$) functions on
$U$. Then, for any $\lambda_0^k\in J^{k}(\bR^q,\bR^m)_{0}$ such
that $(\lambda_0^k,0)\in U$, there exists a uniquely determined
smooth function $F\:U_0\times V_0\to \bR^m$, where $U_0$ is an open
neighborhood of $\lambda^k_0\in J^{k}(\bR^q,\bR^m)_{0}$ and $V_0$
is an open neighborhood of $0\in \bR^q$, such that if $f=(f_1,\ldots,
f_m)$ solves
$$
\partial^\alpha_x f_j=r^\alpha_j(\partial^\beta_x f), \quad\forall
|\alpha|=k+1,\ j=1,\ldots,m, \tag 3.55
$$
near $0\in \bR^q$ and $j^k_0(f)\in U_0$, then
$$
F(j^k_0(f),\cdot)=f.\tag 3.56
$$
\endproclaim

\remark{Remark $3.57$} Observe that we do not claim that
$F(\lambda^k,\cdot)$ solves \thetag{3.55} for any initial value
$\lambda^k$, but only that {\it if} there is a solution with this
initial condition then it coincides with $F(\lambda^k,\cdot)$. The
idea for Proposition 3.54 was suggested to the author by D. Zaitsev.
\endremark
\demo{Proof of Proposition $3.54$} By a standard argument (considering
the derivatives $\partial^\beta_xf$, $|\beta|\leq k$, as new
unknowns), it suffices to prove Proposition 3.54 with $k=1$. Thus, we
may assume that the system \thetag{3.55} is of the form
$$
\partial_{x_j} f_i=r_{ij}(f),\quad i=1,\ldots, m;\ j=1,\ldots, q.\tag 3.58
$$
Fix $\lambda^1_0$ as in the theorem. Write $x=(t,x')\in
\bR\times\bR^{q-1}$ and consider the initial value 
problem for a system of ordinary differential equations
$$
\partial_t f_i(t,0)=r_{i1}(f(t,0))(t,0),\quad f(0,0)=\lambda^1,\tag
3.59
$$
for $\lambda^1$ in some sufficiently small neighborhood of
$\lambda^1_0$. By a classical result (see [CL, Chapter 1.7], Theorem
7.5 and the following remarks), there is
a smooth function $F^1\:U_1\times V_1\to \bR^m$, where $U_1$ is an open
neighborhood of $\lambda^1_0\in J^{1}(\bR^q,\bR^m)_{0}$ and $V_1$
is an open neighborhood of $0\in \bR$, such that $t\mapsto F^1(\lambda^1,t)$ is the
unique solution of \thetag{3.59}. Next, write $x=(x_1,t,x'')\in
\bR\times\bR\times\bR^{q-2}$ and consider for each $x_1\in U_1$ the initial value problem
$$
\partial_t f_i(x_1,t,0)=r_{i1}(f(x_1,t,0))(x_1,t,0),\quad f(x_1,0,0)=F^1(\lambda^1,x_1).\tag
3.60
$$
Again by [CL, Chapter 1.7]
(Theorem 7.5), there is a smooth function $F^2\:U_
2\times V_2\to \bR^m$, where $U_2$ is an open
neighborhood of $\lambda^1_0\in J^{1}(\bR^q,\bR^m)_{0}$ and $V_2$
is an open neighborhood of $(0,0)\in \bR\times \bR$, such that
$t\mapsto F^2(\lambda^1,x_1,t)$ is the 
unique solution of \thetag{3.60}. Proceeding inductively in this way, we obtain
the desired function $F$ after the $q$:th step. The straightforward
details are left to the reader. 
We emphasize however that
the function so obtained need not be a solution of the system
\thetag{ 3.58}, but it satisfies $F(j^1_0(f),\cdot)=f$ whenever $f$ is a
solution. This completes the proof of Proposition 3.54.\qed
\enddemo
The proof of Theorem 4 follows by applying Proposition 3.54 to the
system of differential equations provided by Theorem 2. \qed 
\enddemo

\heading 4.1. Concluding remarks\endheading

\subhead 4.1. Abstract CR manifolds\endsubhead In this paper, we have
considered embedded real hypersurfaces as
abstract manifolds with a(n integrable) CR structure. We have
used the fact
that the CR manifolds are embeddable (i.e.\ the CR structure is
integrable) to choose a basis for the
sections of $\bC TM$ that satisfy certain commutation relations (Lemma
1.33). The author felt that the resulting equations \thetag{1.34}
simplified the computations in the proofs to an extent which, by far,
outweighed the loss of generality in assuming that the manifolds are
embeddable. Without the use of the equations \thetag{1.34}, the key
equations \thetag{2.10-13} take the following form
$$
\aligned
& L_{\bar A}\gamma^E_B+\gamma^E_D R^D_{\bar A B}+\eta^E h_{\bar A
B}=\gamma^D_B\overline{\gamma^C_A} \hat R^E_{\bar C D},
\\
& L_{\bar A}\xi+\xi h_{\bar A}=\xi\overline{\gamma^C_A}\hat h_{\bar
C}+\overline{\gamma^C_A} \eta^D \hat h_{\bar C D},
\\
&L_{\bar A}\eta^C+\eta^C h _{\bar A}+\gamma^C_D R^D_{\bar A}=\xi
\overline{\gamma^E_A} \hat R^C_{\bar E}+\overline{\gamma^E_A}
\eta^F\hat R^C_{\bar E F},
\\
&T\gamma^C_A-L_A\eta^C+\gamma^C_B R^B_A-\eta^C \overline{ h_{\bar
A}}=\xi\gamma^B_A\hat R^C_B+\gamma^B_A \eta^D\hat
R^C_{DB}+\gamma^B_A\overline{\eta^E}\hat R^C_{\bar E
B}.\endaligned\tag 4.1.1 $$ Analogous reflection formulae to those
in Theorem 2.4, as well as analogues of the crucial Lemmas 3.1,
3.6, and Proposition 3.18, can be established (with considerably
more work than above). The author is confident that a proof of
Theorem 2 for abstract CR manifolds (of hypersurface type) $M$ and
$M'$ of the same dimension can be produced from these ingredients,
but he has not had the patience to check the details.

\subhead 4.2. Infinitesimal CR automorphisms\endsubhead It is clear
from the proof of Theorem 3 above that in order for the estimate
\thetag{0.5} to hold, it suffices that $M$ is minimal at $p_0$ and
that there exists an open subset $U\subset M$ with $p_0$ in its
boundary such that $M$ is finitely nondegenerate on $U$. The latter
holds, in particular, if $M$ is holomorphically nondegenerate at
$p_0$ (and, in the real-analytic case, only if), but may hold, in the
case of merely smooth manifolds, even if $M$ is
holomorphically degenerate at $p_0$ (see [BER3, Example 11.7.29]). On
the other hand, as is mentioned in the introduction, if there exists a
vector field $Y$ of the form \thetag{0.4}, where the restrictions of the
$a_j$ to $M$ are
CR functions, which is tangent to $M$ near $p_0$, then
$\dim_\bR\aut(M,p_0)=\infty$. Let us call the restriction to $M$ of such a
vector field $Y$ a {\it CR holomorphic} vector field. Thus, one is led to the following
question.
{\it Suppose $M$ is not finitely nondegenerate at any point in an
open neighboorhood of $p_0$. Does there then exist a CR holomorphic vector
field on $M$ near $p_0$?}
The author does not know the answer to this
question in general, but it seems to be related to the range of the tangential
Cauchy-Riemann operator $\bar\partial_b$ (see e.g.\ [B] for the
definition). We conclude
this paper by briefly outlining this connection.

First, observe that a vector field $Y$ is CR holomorphic if and only if
$Y$ is a section of $\bar\Cal V$ and $[\bar L,Y]$ is a CR vector field
(i.e.\ a section of $\Cal V$) for every CR vector field $\bar L$. Now,
suppose that there is an open set $U\subset M$ in which $M$ is not
finitely nondegenerate at any point. We claim that there exists a (non-vanishing)
CR holomorphic vector field $Y$ near $p\in U$ if (i) $\dim_\bC E_{N-1}(q)$ (which is
$<N$ for $q$ in $U$ by assumption) is maximal at $q=p$, and (ii)
$\bar\partial_b u=f$ is solvable at $p$ for every $(0,1)$-form $f$ with
$\bar \partial_bf=0$. For simplicity, we shall indicate the proof of
this only in the
case $\dim_\bR E_{N-1}(p)=N-1$. We choose a smooth nonvanishing
section $X$ of $\bar \Cal V$ near $p$ such that $X(q)\in E_k(q)^\perp$
for all $k$ and all $q$ in an open neighborhood of $p$. (This can be
done by assumption (i) above.)  We
denote by $L_{\bar 1},\ldots, L_{\bar n}$ a basis of the CR vector
fields on $M$ near $p$, where $L_n:={\overline L_{\bar n}}=X$. We
choose a tranversal vector field $T$, as in \S1, and denote by
$\theta,\theta^A,\bar \theta^{\bar A}$, with notation and conventions
as in \S1, a dual basis of $T,L_A,L_{\bar A}$. The fact that $L_n$ is
valued in $E_k^\perp$, for every $k$, implies that
$[L_{\bar A},L_n]=b_{\bar A}L_n$ modulo the CR vector fields. We shall
look for a CR holomorphic vector field $Y$ of the form $uL_n$, where
$u$ is a function to be determined. It is easy to see that $[L_{\bar
A},Y]$ is a CR vector
field if and only if
$$
L_{\bar A}u+ub_{\bar A}=0\tag 4.2.1
$$
and, hence, $Y$ is CR holomorphic if and only if \thetag{4.2.1} is
satisfied for every $A\in\{1,\ldots,n\}$. If we could solve
$$
\bar \partial_b v=f,\tag 4.2.2
$$
where $f=b_{\bar A}\theta^{\bar A}$, then $u=e^{-v}$ would solve
\thetag{4.2.1}. The $(0,1)$-form $f$ coincides with the form
$L_n\lrcorner\bar\partial_b\theta^n$, as the reader can verify. From
this observation, one can check that $f$ satisfies the necessary
compatibility condition for solvability,
$$
\bar\partial_b f=\bar \partial_b(L_n\lrcorner\bar\partial_b\theta^n)=0.\tag 4.2.3
$$
Hence, if the tangential Cauchy-Riemann complex is solvable at level
$(0,1)$ at $p$, then we can solve \thetag{4.2.2} and obtain a CR
holomorphic vector field $Y=uL_n$ near $p$. However, the author knows of no
results on solvability which apply in this situation (unless, of
course, $M$ is real-analytic).

\enddocument
\end